\documentclass[12pt]{amsart}
\usepackage[left=1in,top=1in,right=1in]{geometry}

\usepackage{
color,
amsmath,
amsfonts,
latexsym, 
amssymb,
amsthm,
hyperref,
ulem,
bbm,
graphicx
}

\usepackage[non-sorted-cites]{amsrefs}

\usepackage{enumerate}

\newcommand{\e}{\varepsilon}
\newcommand{\ii}{\mathrm{i}}

\newcommand{\rd}{{\rm d}}
\newcommand{\m}{\text{\rm meas}}

\renewcommand{\P}{\mathbb{P}}

\newcommand{\E}{\mathbb{E}}
\newcommand{\re}{{\rm Re}}

\DeclareMathOperator{\OO}{O}
\DeclareMathOperator{\oo}{o}

\theoremstyle{plain} 
\newtheorem{theorem}{Theorem}[section]
\newtheorem{lemma}[theorem]{Lemma}

\newtheorem{proposition}[theorem]{Proposition}

\title[Maximum of the Riemann zeta function on a short interval]{Maximum of the Riemann zeta function on a short interval of the critical line}

\author[L.-P. Arguin]{L.-P. Arguin}
\address{Department of Mathematics, Baruch College and Graduate Center, City University of New York, USA}
\email{louis-pierre.arguin@baruch.cuny.edu}

\author[D. Belius]{D. Belius}
\address{Institute of Mathematics, University of Zurich, Switzerland}
\email{david.belius@cantab.net}

\author[P. Bourgade]{P. Bourgade}
\address{Courant Institute, New York University, USA}
\email{bourgade@cims.nyu.edu}

\author[M. Radziwi\l\l]{M. Radziwi\l\l}
\address{Department of Mathematics and Statistics, McGill University, Canada}
\email{maksym.radziwill@mcgill.ca}

\author[K. Soundararajan]{K. Soundararajan}
\address{Department of Mathematics, Stanford University, USA}
\email{ ksound@stanford.edu}

\date{December 14, 2017}

\begin{document}

\begin{abstract}
We prove the leading order of a conjecture by Fyodorov, Hiary and Keating, about the maximum of the Riemann zeta function on random intervals along the critical line. More precisely, as $T \rightarrow \infty$ for a set of $t \in [T, 2T]$ of measure $(1 - \oo(1)) T$, we have 
$$
\max_{|t-u|\leq 1}\log\left|\zeta\left(\tfrac{1}{2}+\ii u\right)\right|=(1+\oo(1))\log\log T .
$$
\end{abstract}

\maketitle


\section{Introduction}

\subsection{Maximum of the Riemann $\zeta$ function on large and short intervals.}\ An important problem in number theory 
concerns the maximum size of the Riemann zeta function on the critical line.   The fundamental Lindel{\" o}f hypothesis \cite{Lin08}
asserts that for any $\varepsilon >0$ and as $|t|\to \infty$ one has $|\zeta(\tfrac{1}{2}+\ii t)|=\OO(|t|^\e)$. 
Among the many arithmetic consequences of the Lindel\"of hypothesis we highlight the existence of primes in all intervals $[x, x + x^{1/2 + \varepsilon}]$ for all $x$ large enough, and in almost all intervals of the form $[x, x + x^{\varepsilon}]$. The current best bound towards the Lindel\"of hypothesis states that $|\zeta(\tfrac 12 + it)| \ll |t|^{13/84 + \varepsilon}$ (see \cite{Bou16}). Chapter XIII of \cite{Tit86} gives a more comprehensive account of the literature surrounding the Lindel{\" o}f hypothesis. 

In \cite{Lit24}, Littlewood showed that a stronger form of the Lindel{\" o}f hypothesis follows from the Riemann hypothesis: namely, for 
some positive constant $C>0$, and for  all large $|t|$ 
\begin{equation}\label{eqn:MaxRH}
|\zeta\left(\tfrac{1}{2}+\ii t\right)|=\OO\left(\exp\left(C\frac{\log |t|}{\log\log |t|}\right)\right). 
\end{equation}
  While the value of the constant $C$ has been reduced over the years \cite{RamSan93, Sou09, ChaSou11}, with \cite{ChaSou11} 
establishing that any $C >(\log 2)/2$ is permissible, Littlewood's  bound remains essentially the best that is known.

There has been more progress on lower bounds for the maximal size of the zeta function. The first result is due to Titchmarsh (see Theorem 8.12 of \cite{Tit86}), 
who proved that for any $\alpha < \tfrac 12$, and large enough $T$, 
$$
\max_{t\in[0,T]}|\zeta(1/2+\ii t)|\geq \exp((\log T)^{\alpha}). 
$$ 
This result was improved to
$$
\max_{t\in[0,T]}|\zeta(1/2+\ii t)|\geq \exp\Big(c\frac{\sqrt{\log T}}{\sqrt{\log\log T}}\Big)
$$ 
in  \cite{Mon77} under the Riemann hypothesis, and then made unconditional with improved constant $c$ in \cite{BalRam77} and \cite{Sou08}. 
A breakthrough was achieved in recent work of Bondarenko and Seip  \cite{BonSei15} who showed that for any $c<1/\sqrt{2}$,
\begin{equation}\label{eqn:lowerb}
\max_{t\in[0,T]}\left|\zeta\left(\tfrac{1}{2}+\ii t\right)\right|\geq \exp\Big(c\,\frac{\sqrt{\log T\log\log\log T}}{\sqrt{\log\log T}}\Big).
\end{equation}

There is a gulf between the known conditional upper bound (\ref{eqn:MaxRH}) and the unconditional $\Omega$-result (\ref{eqn:lowerb}), and the 
asymptotics of the maximal order remains unclear, and a matter of dispute.  In \cite{FarGonHug07}, Farmer, Gonek and Hughes have conjectured that 
$$
\max_{t\in[0,T]}\log |\zeta(1/2+\ii t)|\sim \sqrt{\tfrac{1}{2}\log T\log\log T}, 
$$
but at the end of their paper they also point to dissenting views, advocating that \eqref{eqn:MaxRH} is closer to the maximal size. Extensive numerical computations have been recently carried out in \cite{BoGh}, however the data regarding extreme values remains inconclusive.

Motivated by the goal of understanding the maximum order of $|\zeta(\frac 12+it)|$, Fyodorov, Hiary, and Keating \cite{FyoHiaKea12,FyoKea14} 
proposed the study of the maximum size of the zeta function in randomly chosen intervals (on the critical line) of constant length.
They obtained a precise conjecture (supported by numerical data) for the distribution of this maximum over short intervals.
Namely, if $t$ is chosen uniformly from $ [T,2T]$, then
\begin{equation}\label{eqn:FyoHiaKea}
\max_{|t-u|\leq 1}\log \left|\zeta\left(\tfrac{1}{2}+\ii u\right)\right|=\log\log T-\tfrac{3}{4}\log\log\log T+X_T,
\end{equation}
where the random variable $X_T$ converges weakly, as $T\to\infty$, to an explicitly given distribution.   
Here, for convenience, we have stated their conjecture for random intervals of length $2$, but a similar conjecture 
could be made for intervals of any constant length.    The main result  of this paper is a proof of the leading order asymptotics in (\ref{eqn:FyoHiaKea}).

\begin{theorem}\label{thm:main}
For any $\e>0$, as $T\to \infty$ we have
$$
\frac{1}{T}\m \Big\{ T\leq t\leq 2T: (1-\e)\log\log T<\max_{|t-u|\leq 1}\log\left|\zeta\left(\tfrac{1}{2}+\ii u\right)\right|<(1+\e)\log \log T\Big\} \to 1.
$$
\end{theorem}

While completing this work, we learned that Theorem \ref{thm:main} (as well as the analogue for ${\rm Im} \log\zeta$) was independently proved by Najnudel \cite{Naj16} under the assumption of the Riemann hypothesis. 
It  would be interesting to establish the result for ${\rm Im} \log \zeta$ unconditionally, perhaps by a modification of the approach given here.

\subsection{Extrema of log-correlated fields.}\ 
Fyodorov, Hiary and Keating's conjecture was motivated by a connection with random matrices. This analogy has been the subject  
of many investigations, beginning with  Montgomery's pair correlation conjecture \cite{Mon73}, and leading more recently to the Keating--Snaith conjectures about the moments of the Riemann zeta function \cite{KeaSna00}.   While the pair correlation conjecture 
examines this analogy on the ``microscopic" scale of the average spacing between consecutive zeros (which is $1/\log T$ at 
height $T$),  the prediction (\ref{eqn:FyoHiaKea}) relies on the analogy at a larger ``mesoscopic" scale (intermediate between $1/\log T$, and the ``macroscopic" scale of size $1$).    

To give a sense of this, we recall the fundamental result of Selberg \cite{Sel46} that  if $t$ is sampled uniformly at random from $[T,2T]$ then $\log|\zeta(\tfrac 12+ \ii t)|$  is normally distributed with mean $\sim 0$, and  variance $\sim \frac{1}{2}\log\log T$. His central limit theorem has been extended to study the correlation between values of the zeta function at nearby points
in \cite{Bou10}: for example, if $t$ is uniform on $[T,2T]$ and $0<h<1$, then the covariance between $\log|\zeta(\frac{1}{2}+\ii t)|$ and $\log|\zeta(\frac{1}{2}+\ii (t+h))|$ is 
\begin{equation}\label{eqn:convarianceform}
\tfrac{1}{2}\log\min \left(h^{-1} ,{\log T}\right). 
\end{equation}
Here the comparison of $h^{-1}$ and $\log T$ is natural since $1/\log T$ is (as mentioned above) the scale of the typical 
spacing between zeros of $\zeta(s)$.

A parallel story holds for the logarithm of the characteristic polynomial of $N\times N$ Haar-distributed unitary matrices, 
$\log |{\rm P}_N(z)|$.  On the unit circle $|z|=1$, the distribution of  $\log |{\rm P}_N(z)|$ is asymptotically Gaussian with mean $0$ 
and variance $\sim \frac{1}{2}\log N$ \cite{KeaSna00}.  Moreover, for two points $z_1$ and $z_2$ on the unit circle within distance 
$|z_1-z_2|=h$, the covariance between $\log |{\rm P}_N(z_1)|$ and $\log |{\rm P}_N(z_2)|$ is roughly $\frac{1}{2}\log\min \left(h^{-1},{N}\right)$, analogously to (\ref{eqn:convarianceform}) (see \cite{Bou10}). Fyodorov, Hiary and Keating gave a very precise conjecture for the maximum of $\{\log|{\rm P}_N(z)|,|z|=1\}$ by relying on the replica method, and techniques from statistical mechanics predicting extreme values in disordered systems \cite{FyoBou08,FyoLedRos09,FyoLedRos12}. Assuming that the structure of the logarithmic covariance governs the distribution of the extreme values of $\log |{\rm P}_N(z)|$, they were led to conjecture the
asymptotics (\ref{eqn:FyoHiaKea}).

The above Fyodorov-Hiary-Keating picture of extreme value theory has recently been proved in a variety of cases. For a probabilistic model of the Riemann zeta function the leading order of the maximum on short intervals was obtained in \cite{Har13}, and the second order in \cite{ArgBelHar15}. For the characteristic polynomial of random unitary matrices, the asymptotics of the maximum at first order \cite{ArgBelBou15} and then second order \cite{PaqZei17} are known, together with tightness of the third order \cite{ChaMadNaj16} in the more general context of circular beta ensembles. In the context of Hermitian invariant ensembles, the first order of the maximum of the characteristic polynomial was proved in \cite{LamPaq16} and precise conjectures can be found in \cite{FyoSim2015}.
Theorem \ref{thm:main} and its conditional analogue in \cite{Naj16} are the first results about the maxima of $\zeta$ itself, with the only source of randomness being the choice of the interval.  In connection with the prediction from \cite{FyoHiaKea12,FyoKea14} that $\log|\zeta|$ behaves like a real log-correlated random field, we note that \cite{SakWeb16} recently proved that $\zeta$ converges to a complex Gaussian multiplicative chaos.

To summarize this discussion of related work, we note that our work builds on, and adds to, the efforts to develop extreme value theory of correlated systems. Such statistics are expected to lie on the same universality class for any covariance of type (\ref{eqn:convarianceform}). This class includes the two-dimensional Gaussian free field, branching random walks, cover times of random walks, Gaussian multiplicative chaos, random matrices and the Riemann zeta function. We do not give here a list of the many rigorous works on this topic in recent years, pointing instead to \cite{Arg16, Kis15} and the references therein.

\subsection{About the proof.}\   Theorem \ref{thm:main} asserts two statements: first an {\sl upper bound} that 
for typical $t \in [T,2T]$ one has $\max_{|t-u|\le 1} \log |\zeta(\frac 12 + \ii u)| \le (1+\e) \log \log T$, and second 
a {\sl lower bound} that this maximum is also typically $\ge (1-\e)\log \log T$.  
The upper bound in Theorem \ref{thm:main} admits a short proof based 
on  a Sobolev type inequality and classical second moment estimates for $\zeta(s)$ and $\zeta'(s)$.  This argument 
is given in section \ref{sec:upper}, and indeed in Proposition \ref{prop: upper bound} we 
establish the stronger assertion that for any function $V=V(T)$ tending to infinity with $T$ we have 
$$
\frac{1}{T}\m\Big\{\max_{|t-u|\leq 1}\log\left|\zeta\left(\tfrac{1}{2}+\ii u\right)\right|< \log (V\log T) \Big\}\to 1.
$$
This result is also obtained unconditionally in \cite{Naj16}, by a different argument.  

The lower bound in  Theorem \ref{thm:main} requires substantially more work, and forms the bulk of the paper.  In Section 
3, we reduce the proof of Theorem \ref{thm:main} to two propositions.  The first step, Proposition \ref{prop: lower bound reduction}, 
 transforms the problem to the study of Dirichlet polynomials supported on the primes below $X = \exp((\log T)^{1-\kappa})$ 
 for a suitable $\kappa =\kappa(\varepsilon)>0$.  This reduction step, carried out in Section 4, builds upon ideas from \cite{RadSou15}, which gave an  alternative approach to Selberg's central limit theorem for $\log |\zeta(\frac 12+\ii t)|$.
 The second step, Proposition \ref{prop: lower bound dirichlet}, applies techniques from the theory of branching random walks to establish lower bounds for the Dirichlet polynomials over primes, adapting the approach of \cite{ArgBelBou15, ArgBelHar15}.
 This argument is presented in Section 5.   There is some scope to refine our results by letting the parameter $\kappa$ tend to $0$ (or 
 equivalently the parameter $K$ that will arise later to tend to infinity), but we have not attempted to carry this out.  

 In broad strokes, the proof of Proposition \ref{prop: lower bound reduction} splits into three steps.  
 First we show (Lemma \ref{lem: off axis}) that a large value of $\zeta(s)$ slightly to the right of the critical line (that is, on the 
 line Re$(s) =\sigma_0$ for a suitable $\sigma_0 >\frac 12$)  typically 
 propagates to give a large value on the critical line.  
 In the second step, we construct a finite Dirichlet polynomial $M(s)$ such that for most $t \in [T, 2T]$ and all $u$ with $|t - u| \leq 1$ one has $\zeta(\sigma_0 + \ii u)M(\sigma_0 + \ii u) \approx 1$, with $\sigma_0$ being taken slightly to the right of the half-line 
 (Lemmas \ref{lem4.2} and \ref{lem: cutoff}).  Note that such a construction is not possible if $\sigma_0 = \tfrac 12$ because of the preponderance of zeros of $\zeta(s)$ on the line $\sigma = \tfrac 12$. We call such an $M(s)$ a mollifier. Our mollifier $M(s)$ is constructed in a specific way that allows us in our third step to show that for almost all $t \in [T, 2T]$ we have $M(\sigma_0 + \ii u) \approx \prod_{p \leq X} (1 - p^{-\sigma_0 - \ii u})$ for all $|u - t| \leq 1$, with $X$ substantially smaller than $T$. 
 Assembling together the three steps 
shows that for almost all $t \in [T, 2T]$ a large value of $ \max_{|t - u| \leq 1} \text{Re} \sum_{p\le X} p^{-\sigma_0 -\ii u}$ leads to a large value of $\max_{|t - u| \leq 1} \log |\zeta(\tfrac 12+\ii u)|$.  

We now describe the ideas behind the proof of Proposition \ref{prop: lower bound dirichlet}, where the 
goal is to show that for almost all $t\in [T,2T]$ we have $\max_{|t - u| \leq 1} \text{Re} \sum_{p \leq X} p^{-\sigma_0 - \ii u} >
 (1 - \varepsilon) \log\log T$.   The sketch below is a simplified account of the argument in Section 5, and the reader should be aware 
 of minor discrepancies in notation.  Here $X= \exp((\log T)^{1-\frac 1K})$ for a fixed large integer $K = K(\e)$, and we split the 
 interval $[2,X]$ into $K-1$ disjoint intervals $J_j$ (with $0\le j\le K-2$) setting 
 $J_j = (\exp((\log T)^{\frac jK}, \exp((\log T)^{\frac{j+1}{K}}]$.  Correspondingly, we decompose  $\text{Re} \sum_{p \leq X} p^{-\sigma_0 - \ii u}$ as $\sum_{j=0}^{K-2} P_j(u)$, where the Dirichlet polynomial $P_j$ includes the primes from the interval $J_j$.  
 The interval $J_j$ have been chosen so that for a random $t$ uniformly distributed in $[T, 2T]$, 
\begin{itemize}
\item for $0 \leq j \leq K - 3$, the terms $P_j(t)$ have comparable variance, precisely $\text{var}(P_j(t)) = \frac{1}{2K} (1 + \oo(1)) \log\log T$. 
\item if $j\neq k$ then $P_j(t+\tau)$ and $P_k(t+\tau')$ are asymptotically independent for all fixed $\tau, \tau' \in [0, 1]$.
  \item  for every $j$ and fixed $\tau, \tau' \in [0,1]$, 
\begin{equation}
\label{eqn: split}
{\rm cov}(P_j({t+\tau}),P_j({t+\tau'}))\sim
\begin{cases} 
\frac{1}{2K} \log\log T& {\rm if} -\log |{\tau-\tau'}|\ge \frac{j+1}{K} \log \log T\\
\oo(\log\log T) & {\rm if} -\log |{\tau-\tau'}|\le   \frac{j}{K} \log \log T. 
\end{cases}
\end{equation}
\end{itemize}

 The terms $P_{K-2}(t)$ (which has a slightly different variance from the other terms) and $P_0(t)$ (which correlates 
along fairly long intervals) are special, and it is convenient to discard them. 
This is already anticipated in the statement of Proposition \ref{prop: lower bound reduction}. The next step is to show that for almost all $t \in [T,2T]$ there exists $u$ with $|u-t|\le \frac 14$ and such that $P_j(u) \ge \frac{1-\e}{K} \log \log T$ for all $1\le j\le K-3$.  
 
The Dirichlet polynomials $P_j(t)$ typically do not vary much along intervals of length $1/\log T$, and so one must show 
that for almost all $t \in [T,2T]$ there exists $0\le k < \log T$ with $P_j(t+k/\log T) \ge \frac{1-\e}{K} \log \log T$ for all $1\le j\le K-3$.  
Letting $\mathcal{T}(k / \log T)$ denote the event 
``$P_j(t + k / \log T) > \frac{1-\e}{K} \log\log T$ holds for all $1 \leq j \leq K - 3$,'' an application of the Cauchy-Schwarz inequality gives 
$$
\mathbb{P} \Big ( \bigcup_{0 \leq k < \log T} \mathcal{T}(k / \log T) \Big ) \geq \frac{\Big ( \sum_{0 \leq k < \log T} \mathbb{P}(\mathcal{T}(k / \log T)) \Big )^2 }{ \sum_{0 \leq k, \ell < \log T} \mathbb{P}(\mathcal{T}(k / \log T) \cap \mathcal{T}(\ell / \log T))}.
$$

To evaluate the probabilities arising above, we perform a precise analysis 
in the large deviation regime of the joint distributions of $P_j(t+k/\log T)$ and 
$P_j(t+\ell/\log T)$.  The analysis shows that this joint distribution matches that of 
Gaussian random variables with the covariance structure laid out in \eqref{eqn: split}.
If $k$ and $\ell$ are such that $|k-\ell| > (\log T)^{1-\frac{1}{2K}}$, then for all $1\le j \le K-3$ 
the Dirichlet polynomials $P_j(t+k/\log T)$ and $P_j(k+\ell/\log T)$ behave independently, so that (see Proposition \ref{prop: prob decouple}) 
$$ 
\P ({\mathcal T}(k/\log T) \cap \mathcal{T}(\ell/\log T)) \sim \P({\mathcal T}(k/\log T)) \P (\mathcal{T}(\ell/\log T)).    
$$
Therefore, 
\begin{equation*} \label{eqn:first}
\sum_{\substack{0 \leq k , \ell \leq \log T \\ |k-\ell| > (\log T)^{1-\frac 1{2K}} }} \mathbb{P}(\mathcal{T}(k / \log T) \cap \mathcal{T} (\ell / \log T)) \leq (1 + \oo(1)) \Big ( \sum_{0 \leq k < \log T} \mathbb{P}(\mathcal{T}(k / \log T)) \Big )^2.  
\end{equation*}
This case represents the typical situation when $k$ and $\ell$ range from $0$ to $\log T$.  In the atypical case when $k$ and $\ell$ 
are near each other, $P_j(t+k/\log T)$ and $P_j(t+\ell/\log T)$ will correlate for small values of $j$, and behave independently for larger values (see \eqref{eqn: split} and Proposition \ref{prop: two point couple}).   It follows that 
\begin{equation*} \label{eqn:conclude}
\sum_{\substack{0 \leq k, \ell < \log T \\ |k-\ell| \le (\log T)^{1-\frac{1}{2K}} }} \mathbb{P}(\mathcal{T}(k / \log T) \cap \mathcal{T}(\ell / \log T))  = \oo \Big ( \Big ( \sum_{0 \leq k < \log T} \mathbb{P}(\mathcal{T}(k / \log T)) \Big )^2 \Big ), 
\end{equation*}
and the desired Proposition \ref{prop: lower bound dirichlet} follows.

The approximate correlation behavior of the Dirichlet polynomials $P_j(t + k / \log T)$ and $P_j(t + \ell / \log T)$ has an underlying 
tree structure similar to that of a branching random walk.  
Indeed, an accurate model for $P_j(t + k / \log T)$ can be obtained by considering Gaussian random variables ${ \overline{P}_j(k / \log T)}$ indexed by $\log T$ equidistant points $k / \log T$ on $[0,1]$ and with a dependence structure that we now describe. The points $k / \log T$ are identified with the leaves of a rooted tree with $K-1$ generations indexed by $j$, with each vertex in a generation having approximately $(\log T)^{1/K}$ edges.
One places independent and identically distributed copies of a Gaussian random variable $G_j$ with mean $0$ and variance $(1/2K) \log\log T$ at each edge in generation $j$. 
Given $j$, and a leaf $k / \log T$, the random variable $\overline{P}_j(k / \log T)$ is set to be equal to the random variable $G_j$ that is placed on the path from $k /\log T$ to the root of the tree. Thus given a $j$ and two distinct leaves $k / \log T$, $\ell / \log T$
 the random variables $\overline{P}_j(k / \log T)$ and $\overline{P}_j(\ell / \log T)$ are equal if $-\log |(k - \ell) / \log T| > \frac{j + 1}{2 K} \log\log T$ and independent if $-\log |(k - \ell) / \log T| \leq \frac{j}{2K} \log\log T$, similarly to \eqref{eqn: split}. In fact $(\sum_{k=0}^{K-2} \overline{P}(\tau),\tau \in [-1,1])$ serves as a good model of $(\sum_{k=0}^{K-2} P(t+\tau),\tau \in [-1,1])$. This conceptual picture is explained in detail in \cites{ArgBelHar15,ArgBelBou15} and illustrated in the figure below.

\begin{figure}[h]
\begin{center}
\includegraphics[height=5cm]{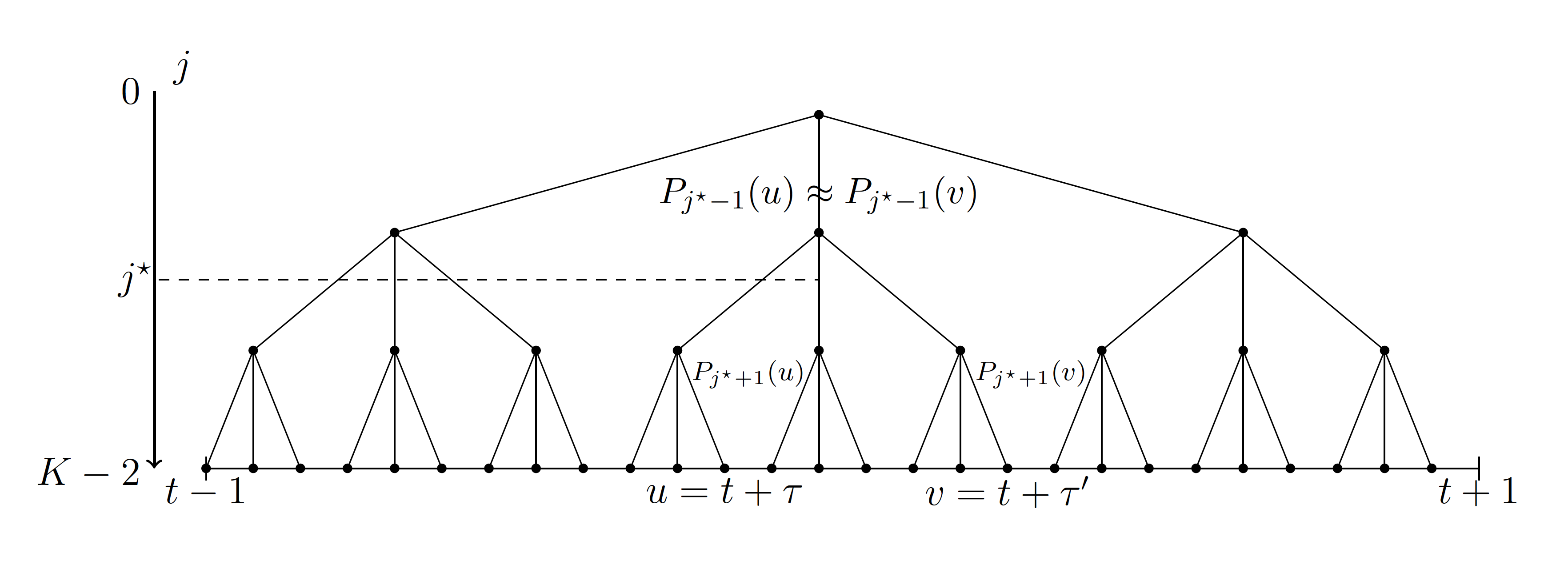}
\caption{Illustration of {the branching random walk approximation} of the Dirichlet polynomials $\sum_{j=0}^{K-2}P_j$. 
For two given {$\tau,\tau'\in[-1,1]$}, the increments {$P_j$} are approximately independent for $j>j^\star$ and {almost identitical} for $j < j^\star$ where $j^\star={K \frac{-\log |\tau-\tau'|}{\log\log T}}$.  {For the model $\sum_{j=0}^{K-2} \overline{P}_j$ this dichotomy holds exactly.}}
\label{fig: tree}
\end{center}
\end{figure}

Finally we remark that the ideas involved in the proof of Proposition 3.2 come from many sources.
The idea of restricting to an initial subset $B$  of $[T, 2T]$ on which an accurate understanding of the distribution of $P_j(t)$ can be obtained comes from \cite{Rad11}.
The identification of an approximate branching random walk structure within the sum $\sum_{p \leq X} p^{-s}$ was used in \cite{ArgBelHar15} to study the extrema of a random model of the zeta function, and in the subsequent works regarding the large values of characteristic polynomials \cite{ArgBelBou15,ChaMadNaj16,LamPaq16,PaqZei17} and of the  zeta function \cite{Naj16}.
The original method for studying the extrema of branching processes which we adapt is due to Bramson \cite{Bra78}.
More precisely, we use Kistler's robust $K$-level coarse graining variant from \cite{Kis15}, as \cite{ArgBelBou15} did for the related random matrix problem.


{\bf Notation.}   For the rest of the paper, { $t$ will denote} a uniform random variable on $[T,2T]$.  Accordingly,  for any event $A_T \subset [T,2T]$ 
and a random variable $X_T: [T,2T] \to {\Bbb C}$ we write 
$$ 
\P(A_T) = \frac 1T \m (A_T), \qquad \text{and} \qquad \E [X_T] = \frac 1T \int_T^{2T} X_T(t) \rd t. 
$$ 
We also use the standard $\OO$ and $\oo$ notations from analytic number theory: thus, $f(T)=\OO(g(T))$ means that $|f(T)|/|g(T)|$ is bounded and $f(T)=\oo(g(T))$ if 
$|f(T)|/|g(T)|\to 0$.    Sometimes it will be convenient to use the notation $f(T) \ll g(T)$, which means the same as $f(T) = \OO(g(T))$.   
We will encounter some arithmetical functions familiar in number theory.  These include: $\omega(n)$ (which 
counts the number of distinct primes dividing $n$), $\Omega(n)$ (which counts with multiplicity the number of primes dividing $n$), 
the von Mangoldt function $\Lambda(n)$ (which equals $\log p$ if $n$ is a power of the prime $p$, and equals $0$ otherwise), and 
 the M{\" o}bius function $\mu(n)$  (which equals $0$ if $n$ is divisible by the square of a prime, and when $n$ is square-free   
 it equals $(-1)^{\omega(n)}$).

{\bf Acknowledgements.} 
The authors thank the referee for useful comments that led to an improvement of the first version of this paper. 
L.-P. A. is supported by NSF CAREER 1653602, NSF grant DMS-1513441, and a Eugene M. Lang Junior Faculty Research Fellowship. D. B. is grateful for the hospitality of the Courant Institute during visits when part of this work was carried out. P. B. is supported by NSF grant DMS-1513587. M. R. is supported by NSERC DG grant, the CRC program and a Sloan Fellowship. K.S.  is partly supported by a grant from the NSF, and a Simons Investigator grant from the Simons Foundation.

\section{Proof of the upper bound}\label{sec:upper}

\noindent The upper bound implicit in our theorem will be a simple consequence of estimates for the 
second moment of the zeta function and its derivative, together with a Sobolev-type inequality.  We begin with 
the Sobolev inequality, which will also be used elsewhere in the proof.   Suppose $f$ (possibly complex valued) is continuously differentiable on $[-1,1]$.  
For any $u \in [-1,1]$, note that 
$$ 
f(u)^2 = \frac{f(1)^2 + f(-1)^2}{2} + \int_{-1}^{u} f'(v) f(v) \rd v - \int_u^1 f^{\prime}(v) f(v) \rd v, 
$$ 
so that using the triangle inequality 
\begin{equation} 
\label{2.1} 
\max_{u\in [-1,1]} |f(u)|^2 \le \frac{|f(1)|^2 + |f(-1)|^2}{2}  + \int_{-1}^{1} |f'(v) f(v)| \rd v. 
\end{equation}

\begin{proposition}
\label{prop: upper bound}
Let $V=V(T)$ be any function that tends to infinity as $T\to \infty$.  Then 
$$
\P\Big(\max_{|t-u|\leq 1}|\zeta(1/2+\ii u)| > V \log T \Big)= \OO (1/V^2) = \oo(1), 
$$
where we recall that $t$ is {sampled} uniformly in the range $[T,2T]$.  
\end{proposition}

\begin{proof}
Chebyshev's inequality implies that
\begin{equation}
\label{eqn: upper cheby}
\P\Big(\max_{|t-u|\leq 1}|\zeta(1/2+\ii u)| > V \log T \Big)
\leq \frac{1}{V^2 (\log T)^{2}} \E\left[\max_{|t-u|\leq 1}|\zeta(1/2+\ii u)|^2\right].
\end{equation}
Applying \eqref{2.1} with $f(v) = \zeta(1/2+\ii t +\ii v)$, we obtain 
$$ 
 \max_{|t-u|\leq 1}|\zeta(1/2+\ii u)|^2  \ll  
  |\zeta(\tfrac 12+\ii (t+1))|^2 + |\zeta(\tfrac 12 + \ii (t-1))|^2 + \int_{-1}^{1} |\zeta^{\prime} (\tfrac 12 + \ii (t+v)) \zeta(\tfrac 12 + \ii (t+v))| 
 \rd v,  
 $$ 
 so that 
 $$ 
 \E\left[\max_{|t-u|\leq 1}|\zeta(1/2+\ii u)|^2\right] \ll \frac{1}{T} \int_{T-1}^{2T+1} \Big( |\zeta(\tfrac 12+ \ii t)|^2 + |\zeta'(\tfrac 12+\ii t) 
 \zeta(\tfrac 12+\ii t)| \Big) \rd t. 
 $$ 
Asymptotics for the second moment of the zeta function and its derivatives are well known (see 
Chapter VII of \cite{Tit86} and, in the case of the derivative, \cite{Con88}), and these imply the bounds 
\begin{equation} 
\label{eqn: moments of zeta}
\int_{T-1}^{2T+1} |\zeta(\tfrac 12+\ii t)|^2 \rd t \ll T \log T, \qquad \text{and} \qquad 
\int_{T-1}^{2T+1} |\zeta^{\prime}(\tfrac 12+\ii t)|^2 \rd t \ll T(\log T)^3. 
\end{equation} 
Using these estimates and Cauchy-Schwarz inequality, we conclude that 
$$
\E\left[\max_{|t-u|\leq 1}|\zeta(1/2+\ii u)|^2\right]\ll (\log T)^2, 
$$
which, in view of \eqref{eqn: upper cheby}, yields the proposition. 
\end{proof}

\section{Plan of the proof of the lower bound}\label{sec:lower}

\noindent The lower bound of Theorem \ref{thm:main} will be proved in two main steps.
First, it is shown that the maximum on a short interval of $\log|\zeta|$ is close 
to the maximum of a Dirichlet polynomial supported on primes slightly to the right of the critical line.  
This is the content of Proposition \ref{prop: lower bound reduction}, whose proof builds upon some ideas 
from \cite{RadSou15}.   Second, a lower bound for the maximum of these Dirichlet polynomials 
on an interval is proved using the robust approach of \cite{Kis15} in Proposition  \ref{prop: lower bound dirichlet}.

The following notation will be used throughout the remainder of the paper.
Motivated by \cite{Kis15}, we will fix a large integer $K=K(\e)$ and divide the primes below
\begin{equation}
\label{eqn: X}
X=\exp( (\log T)^{1-\frac{1}{K}} ) 
\end{equation}
into $K-1$ ranges depending on their size, as follows.  Take $J_0 = [2, \exp( (\log T)^{\frac{1}{K}} )]$, 
and for $1\le j \le K-2$ set 
\begin{equation}
\label{eqn: J}
J_j = (\exp( (\log T)^{\frac{j}{K}} ), \exp( (\log T)^{\frac{j+1}{K}})]. 
\end{equation}
For each $0\le j\le K-2$, we define the Dirichlet polynomial
\begin{equation}
\label{eqn: P}
P_j(u)= \re \sum_{p\in J_j} \frac{1}{p^{\sigma_0+ \ii u}}, 
\end{equation}
where  
\begin{equation}
\label{eqn: sigma}
 \sigma_0 = \frac{1}{2} + \frac{(\log T)^{\frac{3}{2K}}}{\log T}.
 \end{equation}

In the course of the proof, we shall see that if $u$ is chosen uniformly from $[T,2T]$ then $P_j(u)$ {asymptotically} has  a
Gaussian distribution with mean $0$, and variance $\sim \frac 12\sum_{p\in J_j} 1/p^{2\sigma_0}$, see for example Lemma \ref{lem: moments}.   
The prime number 
theorem enables us to evaluate this variance asymptotically, and we record the relevant estimates for future use.  Thus, using 
the prime number theorem (see for example Theorem 6.9 of \cite{MoVa07}) and partial summation it follows that 
for some constant $c>0$, and any $\sigma = \frac 12 +\delta$ with $\delta >0$   
\begin{equation}
\label{eqn: PNT}
\sum_{x \leq p\leq y}\frac{1}{p^{2\sigma}}=\int_x^y \frac{1}{u^{2\sigma}\log u} \rd u+ \OO(e^{-c\sqrt{\log x}}) 
= \log \frac{\log y}{\log x} + \OO (\delta \log y + e^{-c\sqrt{\log x}}).   
\end{equation}
Since $(\sigma_0-1/2)\times\log (\sup J_{K-3}) = (\log T)^{-\frac{1}{2K}}$ it follows that, for all $0\le j\le K-3$, 
\begin{equation} 
\label{3.6} 
\sum_{p \in J_j} \frac{1}{p^{2\sigma_0}} = \frac 1K \log \log T + \OO((\log T)^{-\frac{1}{2K}}), 
\end{equation} 
so that the Dirichlet polynomials $P_j(u)$ all have roughly the same variance.  The last Dirichlet polynomial $P_{K-2}(u)$ 
has a slightly different variance, with the corresponding sum in \eqref{3.6} being roughly  $\frac{1}{2K} \log \log T$.

We are now ready to state the two main propositions from which the lower bound in the theorem  will follow.  
\begin{proposition}
\label{prop: lower bound reduction}  Let $\e>0$ be given, and let $K=K(\e)$ be a suitably 
large integer.  Then  
$$
\begin{aligned}
\P\Big(\max_{|t-u|\leq 1}\log |\zeta(\tfrac 12+\ii u)| &>(1-2\e) \log\log T\Big) \\
&\geq \P\Big(\max_{|t-u|\leq \frac{1}{4}}\sum_{j=1}^{K-3} P_j(u)>(1-\e)\log\log T\Big)+\oo(1).
\end{aligned}
$$
\end{proposition}

Note that in Proposition \ref{prop: lower bound reduction} we omitted the first and last terms, $P_0(u)$ and 
$P_{K-2}(u)$.   The term $P_{K-2}$ is omitted in view of its slightly different variance.
{The very 
small primes occurring in $P_0$ are omitted so that the Dirichlet sums are not too correlated, a fact essential to the analysis in Section 5. }

\begin{proposition}
\label{prop: lower bound dirichlet}  Let $K >3$ be a natural number, and $0< \lambda <1$ be a real number.  Then 
 \begin{equation}\label{eqn:numb1}
\P\Big( \max_{|t-u| \le \frac 14} \Big( \min_{1\le j\le K-3}   P_j(u) \Big) > \frac{\lambda}{K} \log \log T  
\Big)
=1+\oo(1).
\end{equation}
\end{proposition}

\begin{proof}[Proof of Theorem \ref{thm:main}]  If the event of Proposition \ref{prop: lower bound dirichlet} holds,  then 
$$
\max_{|u-t| \le \frac 14} \sum_{j=1}^{K-3}P_j(u)> \lambda\Big(1-\frac{3}{K}\Big)\log\log T.  
$$
Taking $\lambda$ sufficiently close to $1$, and $K$ large enough, the 
  lower bound of the theorem now follows from Proposition \ref{prop: lower bound reduction}.
\end{proof}

Before proceeding to the proofs of the proposition, we record some simple results on 
mean values of Dirichlet polynomials which will be repeatedly used below.  

\begin{lemma} 
\label{lem3.3}  
For any complex numbers $a(n)$ and $b(n)$, and $N\le T$ we have 
$$
\int_{T}^{2T} \Big(\sum_{m\le N} a(m) m^{-\ii t} \Big) \Big( \sum_{n\le N} b(n) n^{\ii t} \Big) \rd t =  T \sum_{n\le N} a(n) b(n) + \OO \Big( N\log N 
\sum_{n\le N} (|a(n)|^2 + |b(n)|^2) \Big).
$$ 
\end{lemma} 
\begin{proof}  Expanding out, and performing the integral, gives  
$$  
 \sum_{m, n\le N} a(m) b(n) \int_T^{2T} \Big( \frac{n}{m} \Big)^{\ii t} \rd t = 
T  \sum_{n\le N} a(n)b(n) + \OO \Big( \sum_{m\neq n \le N} \frac{|a(m)b(n)|}{|\log (m/n)|}  \Big). 
$$ 
Using $|a(m)b(n)| \le |a(m)|^2 +|b(n)|^2$, the remainder term above is 
$$ 
\ll \sum_{m\le N} |a(m)|^2 \sum_{\substack{n\le N \\ n\neq m}} \frac{1}{|\log (m/n)|} + \sum_{n\le N} |b(n)|^2 
\sum_{\substack{ m\le N \\ m\neq n}} \frac{1}{|\log (m/n)|}  \ll N \log N \sum_{n\le N} ( |a(n)|^2+b(n)^2),
$$
proving the lemma. 
\end{proof}

The next two lemmas are also standard (for example, see 
Proposition 3.1 of \cite{Bou10}, or Lemma 3 of \cite{Sou09}), and will be useful in comparing moments of Dirichlet polynomials over the primes 
with the moments of suitable Gaussian distributions.  
\begin{lemma}
\label{lem: moments}  
Let $x\ge 2$ be a real number, and suppose that for primes $p\le x$,  $a(p)$ and $b(p)$ are complex numbers with $|a(p)|$ and $|b(p)|$ 
both at most $1$.   Then for any natural number $k$  we have
$$
\E\Big[
\Big( \frac{1}{2}\sum_{ p \le x} (a(p)p^{-\ii t}+ b(p)p^{\ii t})\Big)^{k}
\Big]
=\partial_z^{k}\Big(\prod_{p\le x} I_0(\sqrt{a(p)b(p)}z)\Big)\Big|_{z=0} + \OO\Big(\frac{x^{2k}}{T}\Big)
$$
where 
$I_0(z)=\sum_{n\geq 0}z^{2n}/(2^{2n}(n!)^2)$ denotes the Bessel function.  
In particular, the expression is $\OO\left(x^{2k}/T\right)$ for odd $k$.
\end{lemma}
\begin{proof}  Given $n$ with prime factorization $n= p_1^{\alpha_1} \cdots p_r^{\alpha_r}$, we set 
$a(n) = \prod_j a(p_j)^{\alpha_j}$, and $b(n) = \prod_j b(p_j)^{\alpha_j}$, with the understanding that $a(n)$ 
and $b(n)$ are $0$ if $n$ has a prime factor larger than $x$.   We also define temporarily the multiplicative function $g$ 
given by $g(p^{\alpha}) =1/\alpha!$ on prime powers $p^{\alpha}$.  With this notation, we may expand (recall $\Omega(n)$ 
counts with multiplicity the number of prime factors of $n$) 
\begin{align*}
\Big( \frac{1}{2}\sum_{ p \le x}& (a(p)p^{-\ii t}+ b(p)p^{\ii t})\Big)^{k}  \\
&= \frac{1}{2^k} \sum_{\ell =0}^{k} \binom{k}{\ell} \Big( \ell!\sum_{\Omega(m) =\ell} a(m)g(m) m^{-\ii t}  \Big) 
\Big( (k-\ell)!\sum_{\Omega(n) = k-\ell} b(n) g(n) n^{\ii t} \Big). 
\end{align*} 
Now we appeal to Lemma \ref{lem3.3} to evaluate the expectation of the above quantity.  The remainder terms 
that arise are 
$$ 
\ll \frac{x^k \log (x^k)}{T} \frac{1}{2^k} \sum_{\ell =0}^{k } \binom{k}{\ell} \Big( \pi(x)^{\ell} + \pi(x)^{k-\ell} \Big) \ll 
\frac{x^{2k}}{T},
$$ 
{ where $\pi(x)\sim x/\log x$ denotes the number of primes below $x$.}

Now let us consider the main terms arising from Lemma \ref{lem3.3}.  These arise from the diagonal terms $m=n$, 
so that $\ell = \Omega(m) = k-\ell = \Omega(n)$.  Thus when $k$ is odd there is no main term, and when $k$ is even, 
we get a main term contribution of 
$$ 
\frac{1}{2^k} k! \sum_{\Omega(n) =k/2} a(n)b(n) g(n)^2.
$$ 
This is $k!$ times the coefficient of $z^k$ in 
$$ 
\sum_{n} a(n) b(n) g(n)^2 (z^2/4)^{\Omega(n)} = \prod_{p\le x} I_0(\sqrt{a(p)b(p)}z),
$$
{since the terms appearing on the left side are multiplicative.}
\end{proof}  

The last lemma gives a useful bound for the $2k$-th moment of Dirichlet polynomials 
supported on primes; it may be deduced by a variant of our argument for the 
previous lemma, or see Lemma 3 of \cite{Sou09}.

\begin{lemma}
\label{lem: bound moments sound}   Let $x\ge 2$ be a real number, and suppose $\sigma \ge \frac12$. 
Let $k$ be a natural number such that $x^k\leq T(\log T)^{-1}$.  Then, for any sequence of complex numbers $a(p)$ 
defined on the primes $p$ below $x$,   
$$
\frac 1T \int_{T}^{2T} \Big | \sum_{p \leq x} \frac{a(p)}{p^{\sigma + \ii t}} \Big |^{2k} \rd t 
\ll  k!  \Big ( \sum_{p \leq x} \frac{|a(p)|^2}{p^{2\sigma}} \Big )^{k}.  
$$
\end{lemma}

\section{Proof of Proposition \ref{prop: lower bound reduction}}

\subsection{Step 1.}   We divide the proof of the proposition into three parts, the first of which bounds the 
maximum of the zeta function over intervals of the critical line in terms of the maximum over intervals 
lying slightly to the right of the critical line.

\begin{lemma}
\label{lem: off axis}  Let $\e>0$ be given, and suppose $\tfrac 12 \leq \sigma \leq \frac{1}{2}+(\log T)^{-1/2 -\e}$. 
Then, for any real number $V\ge 2$,     
$$
\P\Big( \max_{|t-u|\leq 1} |\zeta\left(1/2+\ii  u\right)|>V\Big)
\geq 
\P\Big(\max_{|t-u|\leq \frac{1}{4}} |\zeta(\sigma+\ii u)|>2V\Big)+\oo(1).
$$
\end{lemma}
\begin{proof}  From Theorem 4.11 of \cite{Tit86} we recall that for $\sigma \ge \tfrac 12$ 
\begin{equation} 
\label{4.1} 
\zeta(\sigma+ \ii t) = \sum_{n\le T} \frac{1}{n^{\sigma+\ii t}} + \OO(T^{-\frac 12}). 
\end{equation} 
Using knowledge of the Fourier transform of the function $e^{-|x|}$, we may write 
$$ 
\frac{1}{n^{\sigma -\frac 12}} = \frac 1\pi \int_{-\infty}^{\infty} n^{-iv} \frac{(\sigma -1/2)}{(\sigma -1/2)^2 +v^2} \rd v 
= \frac{1}{\pi} \int_{-T/2}^{T/2} n^{-iv} \frac{(\sigma -1/2)}{(\sigma -1/2)^2 +v^2} \rd v + O(T^{-1}). 
$$ 
Thus we see that 
\begin{equation}
\label{eqn: poisson}
\zeta(\sigma+\ii t)= \frac 1\pi \int_{-T/2}^{T/2} 
\zeta\left(1/2+\ii (t+v) \right)  \frac{\sigma-1/2}{(\sigma-1/2)^2 +v^2}\rd v + O(T^{-\frac 12}).
\end{equation}

Consider $t\in [T,2T]$ such that $\max_{|v|\leq \frac{1}{4}}|\zeta\left(\sigma+\ii (t+v)\right)|>2V$ but 
$\max_{|v| < 1} |\zeta(1/2+\ii (t+v))| \le V$;  we must show that the measure of the set of such points $t$ is $\oo(T)$.  
If $t$ is such a point, then denote by  $v^\star=v^\star(t)$ the $v\in [- \frac{1}{4}, \frac{1}{4}]$ where 
the maximum of $|\zeta(\sigma + \ii (t+v))|$ is attained.  Applying \eqref{eqn: poisson} to the point $\sigma + \ii (t+v^\star)$ 
we obtain 
 $$
2V< |\zeta(\sigma+\ii (t+v^\star))|
\leq \frac{1}{\pi} \int_{-T/2}^{T/2}  |\zeta\left(1/2+\ii (t+v^\star+v)\right)| \frac{(\sigma -1/2)}{(\sigma -1/2)^2 +v^2} \rd v + \OO (T^{-\frac 12}).
$$
Since $|\zeta(1/2+iu)| \le V$ for $|t-u| \le 1$ (by assumption), the portion of the integral above with $|v| \le \frac 34$ is less than
$V$.  Therefore it follows that 
$$ 
V+ \OO(T^{-\frac 12}) \le \frac 1\pi \int_{\frac 34 \le |v| \le \frac{T}{2}} |\zeta\left(1/2+\ii (t+v^\star+v)\right)| \frac{(\sigma -1/2)}{(\sigma -1/2)^2 +v^2} \rd v.
$$  
Using the Cauchy-Schwarz inequality, we deduce that for such $t$, 
$$ 
\Big( \frac{V}{(\sigma -1/2)}\Big)^2 \ll \Big( \int_{\frac 34 \le |v| \le \frac T2} |\zeta(1/2+ \ii (t+v^\star+v))| \frac{\rd v}{v^2} \Big)^2 
\ll \int_{\frac 12 \le |v| \le \frac{T}{2}} |\zeta(1/2+\ii (t+v))|^2 \frac{\rd v}{v^2}. 
$$ 
Therefore, by Chebyshev's inequality, the measure of the set of such points $t\in [T,2T]$  is 
$$ 
\ll \Big( \frac{(\sigma -1/2)}{V}\Big)^2 \int_T^{2T} \int_{\frac 12 \le |v| \le \frac{T}{2}} |\zeta(1/2+\ii (t+v))|^2 \frac{\rd v}{v^2} \rd t 
\ll \Big( \frac{(\sigma -1/2)}{V}\Big)^2 \int_{T/2}^{5T/2} |\zeta(1/2+ \ii t)|^2 \rd t, 
$$ 
which, by \eqref{eqn: moments of zeta} and the assumption on $\sigma$, is 
$$ 
\ll (\sigma -1/2)^2 T \log T = \oo(T). 
$$ 
\end{proof}

\subsection{Step 2.}  The second part of the attack will consist of showing that on the $\sigma_0$ line, one 
can typically invert $\zeta(\sigma_0+it)$ and replace it by a suitable Dirichlet 
polynomial.  
 We define 
\begin{equation}
\label{eqn: M}
M(s)=\sum_{n} \frac{\mu(n)a(n)}{n^{s}},
\end{equation}
where the factor $a(n)$ equals $1$ if all primes factors of $n$ are smaller than $X$  
and $\Omega(n)\leq 100 K\log\log T=:\nu$, and $a(n)=0$ otherwise.    Recall that $\mu$ denotes the M{\" o}bius function, $\Omega(n)$ 
counts the number of prime factors of $n$ (with multiplicity), and $X$ was defined in \eqref{eqn: X}.  
The choice of the Dirichlet polynomial $M$ is 
motivated by work in \cite{RadSou15}, which in turn is motivated by classical ideas on mollifying the zeta function.  
Adapting the proof of Proposition 3 in \cite{RadSou15}, we first establish the following preliminary result. 

\begin{lemma} 
\label{lem4.2}  With the above notation 
$$ 
\int_T^{2T} |\zeta(\sigma_0+\ii t) M(\sigma_0 +\ii t) -1|^2 \rd t = \OO\Big( \frac{T}{(\log T)^{100}}\Big). 
$$ 
\end{lemma}  
\begin{proof}
From its definition, $a(n)=0$ unless $n\le X^{\nu} < T^{\e}$ ($\e>0$ is a fixed arbitrarily small constant), and therefore estimating 
trivially one has $M(\sigma_0+\ii t) \ll T^{\e}$.  Combining this with \eqref{4.1}, we 
see that 
$$ 
\int_T^{2T} \zeta(\sigma_0+\ii t) M(\sigma_0+\ii t) \rd t  
= \int_{T}^{2T} \sum_{n\le T} \frac{1}{n^{\sigma_0+\ii t}} \sum_{m} \frac{\mu(m) a(m)}{m^{\sigma_0+\ii t}} \rd t + \OO (T^{\frac 12+\e}). 
$$ 
Carrying out the integral over $t$, this is 
$$ 
T + \OO\Big( \sum_{\substack{n\le T, m\le X^{\nu} \\ {mn >1 }}} \frac{1}{(mn)^{\sigma_0}} \Big) + \OO(T^{\frac 12+\e}) = T+ \OO(T^{\frac 12+\e}). 
$$ 
Thus, expanding out the square in the desired integral, we see that it equals 
\begin{equation} 
\label{4.4} 
\int_T^{2T} |\zeta(\sigma_0+\ii t)M(\sigma_0 + \ii t)|^2 \rd t  - T + \OO(T^{\frac 12+ \e}). 
\end{equation} 

To estimate the second moment in \eqref{4.4}, we invoke Lemma 4 for \cite{RadSou15}:  
for any $h,k\leq T$ and $1/2<\sigma\leq 1$, we have
\begin{align}
\int_T^{2T}\Big(\frac{h}{k}\Big)^{\ii t}\left|\zeta\left(\sigma+\ii t\right)\right|^2\rd t
&=
\int_T^{2T}\Big(
\zeta(2\sigma)\Big(\frac{(h,k)^2}{hk}\Big)^\sigma
+
\Big(\frac{t}{2\pi}\Big)^{1-2\sigma}\zeta(2-2\sigma)\Big(\frac{(h,k)^2}{hk}\Big)^{1-\sigma}
\Big)\rd t \nonumber\\
&\hskip 1 in +\OO(T^{1-\sigma+\e}\min(h,k)).
\label{eqn:esti}
\end{align}
Using this result, we may write 
\begin{align} 
\label{4.5} 
\int_T^{2T}  |\zeta(\sigma_0+\ii t)M(\sigma_0 + \ii t)|^2 \rd t  &= \sum_{h,k} 
\frac{\mu(h)a(h)\mu(k)a(k)}{h^{\sigma_0} k^{\sigma_0}} \int_T^{2T} \Big(\frac{h}{k}\Big)^{\ii t} |\zeta(\sigma_0+\ii t)|^2 \rd t \nonumber \\
&=  S_1 + S_2 + E, 
\end{align} 
say, with 
\begin{equation} 
\label{4.6} 
S_1= T\zeta(2\sigma_0) \sum_{h, k} \frac{\mu(h)a(h)\mu(k)a(k)}{(hk)^{\sigma_0}} \Big( \frac{(h,k)^2}{hk}\Big)^{\sigma_0}, 
\end{equation} 
\begin{equation} 
\label{4.7} 
S_2 =\zeta(2-2\sigma_0) \Big( \int_T^{2T} \Big( \frac{t}{2\pi}\Big)^{1-2\sigma_0} \rd t \Big)\sum_{h, k} 
\frac{\mu(h) a(h) \mu(k) a(k)}{(hk)^{\sigma_0}} \Big( \frac{(h,k)^2}{hk}\Big)^{1-\sigma_0}, 
\end{equation} 
and 
\begin{equation} 
\label{4.8} 
E = \OO\Big( T^{\frac 12+ \e} \sum_{h, k \le T^{\e} } \frac{1}{(hk)^{\sigma_0}} \min(h,k) \Big)   = \OO(T^{\frac 12+\e}). 
\end{equation} 

Now consider the quantity $S_1$.  Here the sum is over all $h$ and $k$ whose prime factors are below $X$, and with 
$\Omega(h)$ and $\Omega(k)$ below $\nu$.  If we retain the first condition, but drop the second condition, then 
the contribution to $S_1$ would be (upon considering whether a prime $p$ divides neither $h$ nor $k$, or divides exactly one 
of $h$ or $k$, or divides both $h$ and $k$) 
\begin{align}
\label{4.9} 
T\zeta(2\sigma_0) \sum_{\substack{h, k \\ p |hk \implies p\le X}} \frac{\mu(h) \mu(k)}{(hk)^{\sigma_0}}  \Big(\frac{(h,k)^2}{hk}\Big)^{\sigma_0} &= 
T \zeta(2\sigma_0) \prod_{p\le X} \Big( 1- \frac{1}{p^{2\sigma_0}} - \frac{1}{p^{2\sigma_0}} + \frac{1}{p^{2\sigma_0}}\Big) \nonumber \\ 
&=  T \zeta(2\sigma_0) \prod_{p\le X} \Big( 1- \frac{1}{p^{2\sigma_0}} \Big). 
\end{align}
The difference between $S_1$ and \eqref{4.9} comes from the terms with either $\Omega(h)$ or $\Omega(k)$ being larger than 
$\nu$, and these terms give a contribution bounded by (assuming that $\Omega(h)$ is larger than $\nu$) 
\begin{align*} 
&\ll T\zeta(2\sigma_0) \sum_{\substack{h, k \\ \Omega(h)>\nu \\ p| hk\implies p\le X}} 
\frac{|\mu(h)\mu(k)|}{(hk)^{\sigma_0} } \Big(\frac{(h,k)^2}{hk}\Big)^{\sigma_0} \\
&\ll T \zeta(2\sigma_0) e^{-\nu} \sum_{\substack{h, k  \\ p| hk\implies p\le X}}\frac{|\mu(h)\mu(k)|}{(hk)^{\sigma_0} } \Big(\frac{(h,k)^2}{hk}\Big)^{\sigma_0}  e^{\Omega(h)},
\end{align*}
since $e^{\Omega(h) -\nu} \ge 1$ when $\Omega(h) \ge \nu$, and is non-negative for other terms.   The sum over $h$ and $k$ may 
now be expressed as a product over the primes below $X$, yielding 
$$ 
T\zeta(2\sigma_0) e^{-\nu} \prod_{p\le X} \Big( 1+ \frac{e}{p^{2\sigma_0}} + \frac{1}{p^{2\sigma_0}} + \frac{e}{p^{2\sigma_0}}\Big) 
\ll T (\log T) e^{-\nu} \prod_{p\le X} \Big( 1+ \frac{7}{p}\Big) \ll \frac{T}{(\log T)^{100}}. 
$$ 
Thus 
$$ 
S_1 = T\zeta(2\sigma_0) \prod_{p\le X} \Big(1- \frac{1}{p^{2\sigma_0}}\Big) + O\Big( \frac{T}{(\log T)^{100}}\Big) 
= T \prod_{p> X} \Big(1- \frac{1}{p^{2\sigma_0}}\Big)^{-1} + \OO\Big( \frac{T}{(\log T)^{100}}\Big). 
$$ 
Recalling the definitions of $\sigma_0$ and $X$, we find $(\sigma_0-1/2)\log X = (\log T)^{\frac{1}{2K}}$, and so 
$$ 
\sum_{p> X} \log \Big(1-\frac{1}{p^{2\sigma_0}}\Big)^{-1} \ll \sum_{p>X} \frac{1}{p^{2\sigma_0}} 
\ll X^{-(\sigma_0-1/2)} \sum_{p>X} \frac{1}{p^{\sigma_0+1/2}} \ll (\log T)^{-100},
$$ 
which enables us to conclude that $S_1=T + O(T/(\log T)^{100})$.  

Arguing similarly, we see that 
$$ 
S_2 \sim \zeta(2-2\sigma_0) \Big( \int_T^{2T} \Big(\frac{t}{2\pi}\Big)^{1-2\sigma_0}\rd t\Big) \prod_{p\le X} \Big( 1-\frac 2p + \frac{1}{p^{2\sigma_0}}\Big) 
\ll T^{2-2\sigma_0} \log T \ll \frac{T}{(\log T)^{100}}.
$$
Inserting the evaluation of $S_1$ with the estimates for $S_2$ and $E$ into \eqref{4.5}, and then into \eqref{4.4}, we obtain  the lemma.
\end{proof}

Lemma \ref{lem4.2} ensures that for most $t$ one has $\zeta(\sigma_0+\ii t)M(\sigma_0 +\ii t) \approx 1$, and 
we next refine this to ensure that for most $t$ one has $\zeta(\sigma_0+\ii u)M(\sigma_0 +\ii u) \approx 1$ for all 
$u$ with $|u-t|\le 1$.

\begin{lemma}
\label{lem: cutoff}
For any $\e>0$, we have
$$
\P\Big(
\max_{|t-u|\leq 1}\left|
M(\sigma_0+\ii u)\zeta(\sigma_0+\ii u)- 1\right|>\e
\Big)= \oo(1).
$$
\end{lemma}
\begin{proof}  We deduce this from Lemma \ref{lem4.2} and a Sobolev inequality argument.  Note that 
by \eqref{2.1}, we have 
\begin{align*}
\max_{|t-u| \le 1} |\zeta M(\sigma_0+\ii u) -1|^2 
&\ll |\zeta M(\sigma_0 +\ii (t+1))-1|^2 + |\zeta M(\sigma_0 + \ii (t-1)) -1|^2 \\
&+ \int_{t-1}^{t+1} |\zeta M(\sigma_0 +\ii v)-1| | (\zeta^{\prime} M + \zeta M^{\prime})(\sigma_0+\ii v)| \rd v. 
\end{align*}
 Ignoring the end cases $t\in [T,T+1]$ or $t\in [2T-1,2T]$, by Chebyshev's inequality the probability we 
 want to bound is (using the above estimate) 
 $$ 
 \ll \frac 1T + \frac{1}{\varepsilon^2 T} \int_{T}^{2T} \big( |\zeta M(\sigma_0 +\ii (t+1))-1|^2  + |\zeta M(\sigma_0 +\ii t)-1| | (\zeta^{\prime} M + \zeta M^{\prime})(\sigma_0+\ii t)|  \big)\rd t.
 $$
Applying the Cauchy-Schwarz inequality and Lemma \ref{lem4.2} this is 
$$ 
\ll \frac{1}{\varepsilon^2 (\log T)^{100}} +\frac{1}{\varepsilon^2 (\log T)^{50}} \Big( \frac1T \int_{T}^{2T} \big( |\zeta^{\prime} M|^2 +|\zeta M^{\prime}|^2 \big) (\sigma_0 + \ii t) \rd t\Big)^{\frac 12}.
$$ 
We can bound the last integral above by adapting the argument in \cite{RadSou15}, as we did in the proof of Lemma \ref{lem4.2}.  Or, we can finesse the issue by using the Cauchy-Schwarz inequality once again to bound that term by 
$$ 
 \ll \Big(\frac 1T \int_T^{2T} (|\zeta|^4+ |\zeta^{\prime}|^4) (\sigma_0+\ii t) \rd t \Big)^{\frac 14} 
\Big( \frac 1T \int_T^{2T} (|M|^4 +|M^{\prime}|^4)(\sigma_0 +\ii t) \rd t \Big)^{\frac 14}, 
$$ 
and then use the work of Conrey \cite{Con88}\footnote{To be precise, the work there gives an asymptotic for the fourth moment of $\zeta^{\prime}$ on the critical line (the 
fourth moment for $\zeta$ itself is a classical result of Ingham \cite{Ing28}), but this 
implies the same bound on the $\sigma_0$ line as well.}  
to bound the first factor by $\ll (\log T)^2$, and apply Lemma \ref{lem3.3} to bound the second term by $\ll (\log T)^2$.  
This completes the proof, with a lot of room to spare.   
\end{proof} 

\subsection{Step 3.}  The last stage in our proof involves connecting $\log |M(\sigma_0 + \ii t)|$ (for most $t$) with 
(close relatives) of the Dirichlet polynomials over primes $P_j(t)$.   
For $0\le j\le K-2$, define the Dirichlet polynomials
\begin{equation}
\label{eqn: fancy P}
\mathcal P_j(t)=\sum_{n\in J_j} \frac{\Lambda(n)}{n^{\sigma_0 + \ii t} \log n}, \qquad \text{and} \qquad {\widetilde P}_j(t) = 
\sum_{p \in J_j} \frac{1}{p^{\sigma_0 + \ii t}}.
\end{equation}
Note that $P_j(t) $ is simply the real part of ${\widetilde P}_j(t)$, and the difference between ${\mathcal P}_j$ and ${\widetilde P}_j$ is only in the 
prime powers; estimating the contribution of prime cubes and larger powers trivially we 
see that 
\begin{equation} 
\label{4.11} 
Q(t) = \sum_{j=0}^{K-2} ({\mathcal P}_j(t) - {\widetilde P}_j(t)) = \frac 12\sum_{p\le \sqrt{X}} \frac{1}{p^{2\sigma_0 +2\ii t}} + \OO(1). 
\end{equation} 
Our goal is to show that for most $t$ one has $\max_{|t-u| \le 1} |M(\sigma_0 + \ii u) -\exp(-\sum_{j=0}^{K-2} {\mathcal P}_j(u))|
$ is small, and we begin with the following preliminary lemma.  

\begin{lemma} \label{lem4.4} 
With notation as above,
$$
\P\Big( \max_{|t-u| \le 1} |Q(u)| \ge \log \log \log T\Big) = \oo(1), 
$$ 
and 
$$
\P\Big(\max_{|t-u|\le 1} \max_{0\le j\le K-2} |{\widetilde P}_j(u)| \ge {10} K^{-\frac 12} \log \log T \Big) = \oo(1). 
$$ 
\end{lemma} 
\begin{proof}  The Sobolev inequality \eqref{2.1} gives 
$$ 
\max_{|t-u| \le 1} |Q(u)|^2 \ll |Q(t+1)|^2 + |Q(t-1)|^2 + \int_{-1}^{1} |Q(t+v)Q^{\prime}(t+v)| \rd v, 
$$ 
so that, using Chebyshev's inequality and the Cauchy-Schwarz inequality, 
$$ 
(\log \log \log T)^2 \P\Big( \max_{|t-u| \le 1} |Q(u)| \ge \log \log \log T\Big)  \ll \frac 1T+ \E[|Q(t)|^2] + 
\Big( \E [|Q(t)|^2] \E[|Q^{\prime}(t)|^2]\Big)^{\frac 12}.
$$
 A quick calculation with Lemma \ref{lem3.3} shows that $\E[|Q(t)|^2]$ and $\E[|Q^{\prime}(t)|^2]$ are 
 $\OO(1)$, which yields the first assertion of the lemma.  
 
 Let $\ell$ denote a natural number to be chosen later.  Applying \eqref{2.1} to the function ${\widetilde P}_j(t)^{\ell}$, we 
 obtain 
 $$ 
 \max_{|t-u| \le 1} |{\widetilde P}_j(u)|^{2\ell} \ll |{\widetilde P}_j(t-1)|^{2\ell}  + |{\widetilde P}_{j}(t+1)|^{2\ell} + \ell \int_{t-1}^{t+1} |{\widetilde P}_j(v)|^{2\ell -1} |{\widetilde P}_j^{\prime}(v)| \rd v. 
$$ 
Combining this with Chebyshev's inequality and the Cauchy-Schwarz inequality, we may bound 
$\P( \max_{|u-t|\le 1} |{\widetilde P}_j(u)| \ge 10K^{-\frac 12}\log \log T)$ by 
\begin{equation}
\label{eqn: prob lem 4.4}
\ll  \frac 1T + (10K^{-\frac 12}\log \log T)^{-2\ell} \Big( \E [|{\widetilde P}_j(t)|^{2\ell}] + \ell \Big( \E [|{\widetilde P}_j(t)|^{4\ell-2}] \E[|{\widetilde P}_j^{\prime}(t)|^2 ]\Big)^{\frac 12} \Big). 
\end{equation}
Now an application of Lemma \ref{lem3.3} shows that 
$$ 
\E[|{\widetilde P}_j^{\prime}(t)|^2] \ll \sum_{p \in {J_j}} \frac{(\log p)^2}{p^{2\sigma_0}} \ll (\log T)^2, 
$$ 
and an application of Lemma \ref{lem: bound moments sound} gives 
$$ 
\E[|{\widetilde P}_j(t)|^{4\ell -2}] \ll (2\ell-1)! \Big( \sum_{p\in J_j} \frac{1}{p^{2\sigma_0}} \Big)^{2\ell -1} 
\ll  (\ell K^{-1} \log \log T)^{2\ell -1}. 
$$ 
Upon choosing $\ell = [10 \log \log T]$, we conclude from this and \eqref{eqn: prob lem 4.4} that
$$ 
\P\left( \max_{|u-t|\le 1} |{\widetilde P}_j(u)| \ge 10\log \log T\right) \ll (\log T) \Big( \frac{\ell K^{-1} \log \log T}{100K^{-1} (\log \log T)^2}\Big)^{\ell} \ll (\log T)^{-10}.
$$ 
Using a union bound for each $0\le j\le K-2$, we obtain a 
stronger form of the claimed lemma.  
\end{proof}

We are ready to connect $M(\sigma_0+\ii t)$ with $\exp(-\sum_{j=0}^{K-3} {\mathcal P}_j(t))$ for 
most values of $t$.  
\begin{lemma}
\label{lem: M to P}
We have 
$$
\P\Big(
\max_{|t-u|\leq 1}\Big|
M(\sigma_0+\ii u)- \exp\Big(-\sum_{j=0}^{K-2} \mathcal P_j(u)\Big)
\Big|>  (\log T)^{-2}
\Big)= \oo(1)\ .
$$
\end{lemma}
\begin{proof}  Recalling that $\nu = 100K \log \log T$, we define  the truncated exponential
\begin{equation}
\label{eqn: fancy M}
\mathcal M(t)=\sum_{k\leq \nu}\frac{(-1)^k}{k!} \Big(\sum_{j=0}^{K-2} \mathcal P_j(t)\Big)^k.
\end{equation}
By Lemma \ref{lem4.4}, we know that with probability $1 + \oo(1)$ (in $t$) one has 
$$
\max_{|t-u|\le 1} \Big| \sum_{j=0}^{K-2} {\mathcal P}_j(u) \Big| \le 
\max_{|t -u| \le 1} \Big( |Q(u)| + \sum_{j=0}^{K-2} |{\widetilde P}_j(u)| \Big) 
\le 10K \log \log T. 
$$ 
For such a typical $t$, one has 
$$ 
\max_{|u-t| \le 1} \Big| {\mathcal M}(u) - \exp\Big( -\sum_{j=0}^{K-2} {\mathcal P}_j(u) \Big) \Big| 
\le \sum_{k> \nu} \frac{1}{k!} (10K \log \log T)^k \ll (\log T)^{-100}. 
$$ 
Therefore, the lemma would follow once we establish that 
\begin{equation} 
\label{4.13} 
\P\Big(
\max_{|t-u|\leq 1}\left|
M(\sigma_0+\ii u)- \mathcal M(u)
\right|>(\log  T)^{-3}
\Big)= \oo(1)\ .
\end{equation}

The quantities $M(\sigma_0 + \ii u)$ and ${\mathcal M}(u)$ are almost identical, differing only 
in a small number of terms.  More precisely, if we write ${\mathcal M}(u) = \sum_{n} b(n) n^{-\sigma_0 - \ii u}$, 
then one may check that (i) $|b(n)| \le 1$ always, (ii) $b(n)=0$ unless $n\le X^{\nu}$ is composed only of primes below $X$, and 
(iii) $b(n)=\mu(n)a(n)$ unless $\Omega(n) > \nu$, or if there is a prime $p\le X$ such that $p^k |n$ with $p^k>X$.   
Therefore, an application of Lemma \ref{lem3.3} gives 
$$ 
\E[ |M(\sigma_0+\ii t) - {\mathcal M}(t)|^2 ] 
\ll \sum_{\substack{ p| n \implies p\le X \\ \Omega(n) >\nu}} \frac{1}{n} + \Big( \sum_{\substack{ p\le X \\ p^k >X}} \frac{1}{p^k} \Big) 
\Big( \sum_{p|n \implies p\le X} \frac{1}{n} \Big).  
$$ 
The second term above is $\ll (\log X)/\sqrt{X}  \ll (\log T)^{-100}$.  Since $e^{(\Omega(n)-\nu)/2}$ is 
$\ge 1$ when $\Omega(n)>\nu$, and is positive for all other $n$, we may bound the first term above by 
$$ 
e^{-\nu/2} \sum_{p|n \implies p\le X} \frac{e^{\Omega(n)/2}}{n} \ll (\log T)^{-50K} \prod_{p\le X} \Big(1 + \sum_{j=1}^{\infty} 
\frac{e^{j/2}}{p^j} \Big) \ll (\log T)^{-50}. 
$$ 
We conclude that 
$$ 
\E[ |M(\sigma_0+\ii t) - {\mathcal M}(t)|^2 ]  \ll (\log T)^{-50}. 
$$ 
A simple application of Lemma \ref{lem3.3} also shows that $\E[|M^{\prime}(\sigma_0+\ii t)|^2]$ and $\E[|{\mathcal M}^{\prime}(t)|^2]$ 
are $\ll (\log T)^3$.  {The estimate \eqref{4.13} follows as in Lemmas \ref{lem4.4} and \ref{lem: M to P} by a successive application of the Sobolev inequality \eqref{2.1}, Chebyshev's inequality and the Cauchy-Schwarz inequality, proving the lemma.}
\end{proof}

\subsection{Finishing the proof of Proposition \ref{prop: lower bound reduction}}
It is now simply a matter of assembling the results established above.   From Lemma \ref{lem: off axis} we obtain for any $V\ge 2$ 
$$ 
\P( \max_{|t-u| \le 1} |\zeta(\tfrac 12+\ii u)| \ge V) \ge \P(\max_{|t-u|\le \frac 14} |\zeta(\sigma_0 +\ii u)| \ge 2V) + \oo(1).  
$$ 
By Lemma \ref{lem: cutoff} this quantity is 
$$ 
\ge \P( \max_{|t-u| \le \frac 14} |M(\sigma_0 + \ii u)|^{-1} \ge 4V) + \oo (1), 
$$ 
and by Lemma \ref{lem: M to P} the above is 
$$ 
\ge \P \Big(\max_{|t-u| \le \frac 14} \sum_{j=0}^{K-2} \text{Re }{\mathcal P}_j(u) \ge \log (8V) \Big) + \oo(1). 
$$ 
Invoking Lemma \ref{lem4.4}, we may replace Re${\mathcal  P}_j(u)$ by $P_j(u)$ 
with negligible error, and also discard the terms with $j=0$ and $j=K-2$: thus, the quantity above is 
$$ 
\ge \P \Big( \max_{|t-u| \le \frac 14} \sum_{j=1}^{K-3} P_j(u) \ge \log (8V) + \log \log \log T + 20K^{-\frac 12}\log \log T\Big) 
+ \oo(1). 
$$ 
Taking $V= (\log T)^{1-2\e}$, the proposition follows.

\section{Proof of Proposition \ref{prop: lower bound dirichlet}}

\noindent The proof of the proposition is based on large deviation estimates for $P_j(u)$ (defined in \eqref{eqn: P} and  \eqref{eqn: sigma}), see Propositions \ref{prop: two point couple} and \ref{prop: prob decouple}.
In Section \ref{sect: FL transform}, we estimate the Fourier-Laplace transform of $P_j(u)$ in a wide range, using Lemma \ref{lem: moments} to evaluate moments of Dirichlet polynomials.   The large deviation estimates are then derived by inverting the Fourier-Laplace transforms, in Section \ref{sect: LD estimates}.   The proof of Proposition \ref{prop: lower bound dirichlet} is completed in Section \ref{sect: proof Dirichlet}.

\subsection{Fourier-Laplace Transform of Dirichlet Polynomials}
\label{sect: FL transform}
The first step is to show that the moments of sums of  $P_j$'s are very close to Gaussian moments.   

\begin{proposition}
\label{prop: gaussian moments}  For $1\le j\le K-3$ let $\xi_j$ and $\xi_j^{\prime}$ denote complex numbers 
with $|\xi_j|$, $|\xi_j'|\leq (\log T)^{\frac{1}{16K}}$.   Let $\tau$ denote a real number with $|\tau|\le 1$.   If $n \le (\log T)^{\frac 1{2K}}$ 
is odd then 
$$
\E\Big[\Big(\sum_{j=1}^{K-3} \{ \xi_j P_j(t)+\xi_j'P_j(t+\tau) \} \Big)^{n}  \Big]
 =    \OO(\exp(-(\log T)^{\frac{1}{3K}})).  
 $$ 
If $n\leq (\log T)^{\frac{1}{2K}}$ is even, 
\begin{align}
\label{eqn: moment two point}
\E\Big[\Big(\sum_{j=1}^{K-3} \{ \xi_j P_j(t)+\xi_j'P_j(t+\tau) \} \Big)^{n}  \Big]
& =\frac{n!}{2^{n/2}(n/2)!} \ \Big(\sum_{j=1}^{K-3} \{ s_j^2(\xi_j^2+ {\xi'}^2_j) + 2\rho_j(\tau)\xi_j\xi'_j \} \Big)^{n/2} \nonumber \\ 
&\hskip 1 in+  \OO(\exp(-(\log T)^{\frac{1}{3K}})), 
 \end{align}
where
\begin{equation}\label{eqn: sj and rhoj}
s_j^2=\frac{1}{2}\sum_{p\in J_j}p^{-2\sigma_0} \mbox{  and  } \rho_j(\tau)=\frac{1}{2}\sum_{p\in J_j}p^{-2\sigma_0} \cos(\tau \log p).
\end{equation}
\end{proposition}

Ignoring the remainder term, the moments evaluated in \eqref{eqn: moment two point} correspond exactly with what would happen if 
$P_j(t)$ and $P_j(t+\tau)$ were jointly Gaussian with variance $s_j^2$ and covariance $\rho_j(\tau)$, and with 
$P_j(t)$ and $P_j(t+\tau)$ being uncorrelated with $P_k(t)$ and $P_k(t+\tau)$ when $j\neq k$.   
Recall from \eqref{3.6} that the prime number theorem gives (for $1\le j\le K-3$) 
\begin{equation} 
\label{5.01} 
2 s_j^2  = \frac{\log \log T}{K} + \OO ((\log T)^{-\frac{1}{2K}}). 
\end{equation} 
Moreover, by partial summation the prime number theorem also gives (for $1\le j\le K-3$) 
\begin{equation} 
\label{5.02} 
\rho_j(\tau) = \begin{cases} 
\frac{\log \log T}{2K} + \OO(1) &\text{if  } |\tau| \le (\log T)^{-\frac{j+1}{K}} \\ 
\OO( |\tau|^{-1} (\log T)^{-\frac{j}{K}} ) &\text{if  } 1 \ge |\tau| \ge (\log T)^{-\frac{j}{K}}. 
\end{cases}
\end{equation} 
In particular, we see that the polynomials decorrelate for $j\ge 1$  if the distance $\tau$ is large enough.
The term $P_0$ however remains correlated in a large range of $\tau$, and this is the reason for omitting it 
in Proposition \ref{prop: lower bound dirichlet}.  The range $(\log T)^{-\frac jK} \ge |\tau| 
\ge (\log T)^{-\frac{j+1}{K}}$ can also be handled using the prime number theorem, but we do not 
require this, and will just use the trivial bound $-s_j^2 \le \rho_j(\tau) \le s_j^2$ here.

\begin{proof}[Proof of Proposition \ref{prop: gaussian moments}]   Write 
$$ 
\sum_{j=1}^{K-3} \{\xi_j P_j( t)+\xi_j'P_j( t+\tau) \} = \frac{1}{2}\sum_{p} \{a(p)p^{-it}+a^\star(p)p^{it}\}, 
$$ 
where, for primes $p\in J_j$ with $1\le j\le K-3$, we set  
$$
a(p)=(\xi_j+\xi_j'p^{-i\tau})p^{-\sigma_0} 
\qquad \text{and} \qquad 
a^\star(p)=(\xi_j+\xi_j'p^{i\tau})p^{-\sigma_0}, 
$$ 
and put $a(p)=a^\star(p)=0$ for all other $p$.  We now appeal to Lemma \ref{lem: moments} to evaluate the desired $n$-th 
moment.   In the range $n\le (\log T)^{\frac{1}{2K}}$ the error term in Lemma \ref{lem: moments} is easily seen to be $\ll \exp(-(\log T)^{\frac{1}{3K}})$.  When $n$ is odd there is no main term, completing the proof of this case.  

When $n$ is even, the main term from Lemma \ref{lem: moments} arises as the $n$-th derivative (at $z=0$) of 
$$
\prod_{p} I_0(\sqrt{a(p) a^\star(p)}z)
=\prod_{p} \left(1+\frac{a(p)a^\star(p)z^2}{4}+g_p(z)z^4\right)
$$
for $g_p(z)$ an analytic function in a neighborhood of $0$ with $g_p(z)\ll |a(p)a^\star(p)|^2$. 
 Since $a(p)a^\star(p)=\big\{\xi_j^2+{\xi'_j}^2+ 2\xi_j\xi'_j\cos(\tau\log p)\big\}p^{-2\sigma_0}$ for $p\in J_j$,  
 we may expand the product above as 
$$
\prod_{p}  I_0(\sqrt{a(p) a^\star(p)}z)
=\exp\Big(\frac{z^2}{2}\Big(\frac{1}{2}\sum_j \sum_{p \in J_j}p^{-2\sigma_0}\{ \xi_j^2+{\xi'_j}^2 + 2\xi_j\xi'_j\cos(\tau\log p)\}\Big)\Big)F_X(z), 
$$
for $F_X(z)$ a function which is analytic in a neighborhood of $0$, satisfies $F_X(0)=1$, and whose derivatives at $0$ are uniformly bounded by
$$
\sum_{j=1}^{K-3}\sum_{p\in J_j}|a(p)a^\star(p)|^2\ll (\log T)^{\frac{1}{8K}} \sum_{j=1}^{K-3} \sum_{p\in J_j}p^{-2}\ll (\log T)^{\frac{1}{8K}} \exp(-(\log T)^{\frac{1}{K}})\ .
$$
The claim \eqref{eqn: moment two point} follows from Lemma \ref{lem: moments} by taking the $n$-th derivative
(note that the exponential term is exactly the moment generating function of a Gaussian) and noting
that the terms involving a derivative of $F_X(z)$ contribute at most $\ll \exp(-(\log T)^{\frac{1}{3K} } )$.
\end{proof}

We shall use Proposition \ref{prop: gaussian moments} to compute the Fourier-Laplace transform of $P_j(t)$ and $P_{j}(t+\tau)$ 
in wide ranges.  Since these transforms can be dominated by rare extremely large values of $P_j(t)$, it is necessary to 
introduce a cut-off.   With this in mind, we introduce the set
\begin{equation}
\label{eqn: B}
B=\{T\leq t\leq 2T:  |P_j(t)| \leq (\log T)^{\frac{1}{4K}},\ \text{ for  all }  1\leq j\leq K-3\}.
\end{equation}

\begin{lemma}
\label{lem: barrier}  With $B$ as defined in \eqref{eqn: B}, 
\begin{equation}\label{eqn: bound on B complement}
\mathbb P (B^c)\ll \exp\Big(-\frac{(\log T)^{\frac{1}{2K}}}{\log\log T}\Big).
\end{equation}
\end{lemma}
\begin{proof}  By Chebyshev's inequality and Proposition \ref{prop: gaussian moments} we see that for any even 
$n\le (\log T)^{\frac{1}{2K}}$ 
$$ 
\P \Big( |P_j(t)| \ge (\log T)^{\frac{1}{4K}}\Big) \le (\log T)^{-\frac{n}{4K}} \E\Big[ |P_j(t)|^n \Big] \ll (\log T)^{-\frac{n}{4K}}  \Big( \frac{n}{e}s_j^2\Big)^{n/2}. 
$$ 
Taking $n$ to be an even integer approximately $2K (\log T)^{\frac{1}{2K}}/\log \log T$, we see that this probability is 
$\ll \exp(-(\log T)^{\frac{1}{2K}}/\log \log T)$.   The union bound gives
$$
\mathbb P (B^c)\leq \sum_{j=1}^{K-3} \mathbb P( |P_j(t)| > (\log T)^{\frac{1}{4K}}),
$$
 and since $K$ is fixed, the lemma follows. 
\end{proof}

Given a real number $|\tau| \le 1$, let 
$$
B(\tau)=\{T\leq t\leq 2T: P_j(t+\tau)\leq (\log T)^{\frac{1}{4K}},\ \text{ for all  } 1\leq j\leq K-3\}.
$$ 
Thus $B(\tau)$ is essentially a translate of the set $B= B(0)$, and the bound of Lemma \ref{lem: barrier} 
applies to $\P (B(\tau)^c)$ as well.    On $B$ and $B(\tau)$, we can derive precise bounds for the Fourier-Laplace transforms of 
the $P_j$'s for two points.

\begin{proposition}
\label{prop: fourier}   For $1\le j\le K-3$ let $\xi_j$ and $\xi_{j}^{\prime}$ denote complex numbers with
 $|\xi_j|$, $|\xi_j'|\leq (\log T)^{\frac{1}{16K}}$.   Then 
\begin{equation}\label{eqn: one point exp moment}
\E\Big[\exp \Big( \sum_{j=1}^{K-3}\xi_jP_j(t)\Big) \ {\bf 1}_{B}\Big]
= \exp\Big(\frac{1}{2}\sum_{j=1}^{K-3} \xi_j^2s_j^2\Big)   +\OO(\exp(-(\log T)^{\frac{1}{4K}})).
\end{equation}
Further, for any real number $\tau$ with $|\tau| \le 1$ we have 
\begin{align}\label{eqn: two point exp moment}
\E\Big[\exp \Big( \sum_{j=1}^{K-3}&\xi_jP_j(t)+\xi_j'P_j(t+\tau))\Big) \ {\bf 1}_{B\cap B(\tau)}\Big] \nonumber \\
&= \exp\Big(\frac{1}{2}\sum_{j=1}^{K-3}
\{ s_j^2(\xi_j^2+ {\xi'}^2_j) + 2\rho_j(\tau)\xi_j\xi'_j \}
\Big) +\OO(\exp(-(\log T)^{\frac{1}{4K}})).
\end{align}
\end{proposition}
\begin{proof}
We prove the two-point estimate \eqref{eqn: two point exp moment}, the proof of the one-point estimate \eqref{eqn: one point exp moment} is similar (and simpler).  The approach is similar to the proof of Lemma \ref{lem: M to P}, approximating the exponential using 
many terms in the Taylor expansion, and then invoking the Gaussian moments established in Proposition \ref{prop: gaussian moments}. 

We begin with the following simple observation:  if $z$ is a complex number, and $n$ is a natural number $\ge 10(|z|+1)$ then 
\begin{equation} 
\label{5.1} 
\Big| e^z - \sum_{j=0}^{n} \frac{z^j}{j!} \Big| \le \sum_{j=n+1}^{\infty} \frac{|z|^{j}}{j!} \le \frac{|z|^{n}}{n!} \le e^{-n}. 
\end{equation} 
For brevity, write $P_j$ for $P_j(t)$ and $P_j^{\prime}$ for $P_j(t+\tau)$, and similarly put $B' = B(\tau)$.  On the set $B\cap B'$ 
we have 
$$
\Big|\sum_{j=1}^{K-3} (\xi_j P_j +\xi_j^{\prime}P_j^{\prime}) \Big| \le 2 (K-3)  (\log T)^{\frac{1}{4K}+ \frac{1}{16K}} < (\log T)^{\frac{1}{3K}} -1.
$$  
Therefore, using \eqref{5.1}, with $N = 10 (\log T)^{\frac {1}{3K}}$ we obtain 
\begin{equation} 
\label{5.2} 
\E\Big [ \exp \Big( \sum_{j}\xi_jP_j+\xi_j'P_j'\Big) \ {\bf 1}_{B\cap B'}\Big]= \sum_{n\leq N} \frac{1}{n!} \E\Big[\Big( \sum_{j}\xi_jP_j+\xi_j'P_j'\Big)^n 
 \ {\bf 1}_{B\cap B'}\Big] 
+ \OO ( \exp(- (\log T)^{\frac 1{3K}})).
\end{equation} 

Now we show that the moments restricted to $B\cap B'$ appearing in \eqref{5.2} are very nearly the unrestricted moments to 
which we can use Proposition \ref{prop: gaussian moments}.  The Cauchy-Schwarz inequality gives 
$$
\frac{1}{n!} \E\Big[\Big( \sum_{j}\xi_jP_j+\xi_j'P_j'\Big)^n  \ {\bf 1}_{(B\cap B')^c}\Big]
\leq  \frac{1}{n!} \E\Big[\Big( \sum_{j}\xi_jP_j+\xi_j'P_j'\Big)^{2n}\Big]^{1/2} \times (2\mathbb P(B^c))^{1/2}.
$$
Using Proposition  \ref{prop: gaussian moments} (together with the bounds on $|\xi|_j$, $|\xi_j^{\prime}|$ and $s_j^2$)  and Lemma \ref{lem: barrier}, the above 
is 
$$ 
\ll (\log T)^{n} \exp\Big(- \frac{(\log T)^{\frac{1}{2K}}}{2\log \log T}\Big) \ll \exp( -  (\log T)^{\frac{1}{3K}}). 
$$ 
Therefore for $n\le N$ we have 
\begin{equation} 
\label{5.3} 
\frac{1}{n!} \E\Big[\Big( \sum_{j}\xi_jP_j+\xi_j'P_j'\Big)^n  \ {\bf 1}_{B\cap B'}\Big] =\frac{1}{n!} 
 \E\Big[\Big( \sum_{j}\xi_jP_j+\xi_j'P_j'\Big)^n \Big] + \OO ( \exp( -  (\log T)^{\frac{1}{3K}})). 
 \end{equation}
 
 Now we use Proposition  \ref{prop: gaussian moments} to evaluate the unrestricted moments in \eqref{5.3}.  When 
 $n\le N$ is odd, there is no main term, and the quantity in \eqref{5.3} is bounded by $\ll \exp(-(\log T)^{\frac{1}{3K}})$.  
 When $n = 2m \le N$ is even, then Proposition  \ref{prop: gaussian moments} gives 
 $$ 
 \frac{1}{(2m)!} 
 \E\Big[\Big( \sum_{j}\xi_jP_j+\xi_j'P_j'\Big)^{2m} \Big] 
 = \frac{1}{2^m m!} \Big( \sum_{j} \{ s_j^2 (\xi_j^2 +{\xi_j^{\prime}}^2) + 2\rho_j(\tau) \xi_j \xi_j^{\prime} \} \Big)^m + \OO( \exp(-(\log T)^{\frac{1}{3K}})). 
 $$ 
 Inserting this into \eqref{5.3}, and then into \eqref{5.2}, it follows that 
 \begin{align} 
 \label{5.4} 
\E\Big[\exp \Big( \sum_{j}\xi_jP_j+\xi_j'P_j'\Big) \ {\bf 1}_{ B\cap B'}\Big] &= 
\sum_{m\le N/2} \frac{1}{2^m m!} \Big( \sum_{j} \{ s_j^2 (\xi_j^2 +{\xi_j^{\prime}}^2) + 2\rho_j(\tau) \xi_j \xi_j^{\prime} \} \Big)^m  \nonumber \\ 
&\hskip 1 in + \OO( \exp(-(\log T)^{\frac{1}{4K}})).
 \end{align} 
Since $|\xi_j|$ and $|\xi_j'|$ are bounded by $(\log T)^{\frac{1}{16K}}$, an application of \eqref{5.1} shows that the 
above equals 
$$ 
\exp\Big(\frac{1}{2} \sum_j \left\{ s_j^2(\xi_j^2+ {\xi'}^2_j) + 2\rho_j(\tau)\xi_j\xi'_j \right\} \Big) + \OO(\exp(-(\log T)^{\frac{1}{4K}})), 
$$ 
completing the proof.
\end{proof}
 
\subsection{Large Deviation Estimates}
\label{sect: LD estimates}  Proposition \ref{prop: fourier} can be used to get precise large deviation estimates on the variables $P_j$. 
For  $x_j$ (with $1\le j \le K-3$)  to be fixed later, and $\tau$ a real number with $|\tau|\le 1$, define the events
\begin{equation}
\label{eqn: A}
A(\tau)=\{T\leq t\leq 2T : P_j(t+\tau)> x_j,\ \text{for all }  1\leq j\leq K-3\}.  
\end{equation}
We will abbreviate $A(0)$ as $A$, and note that (away from a bounded distance of the end points $T$ and $2T$) the set $A(\tau)$ is just a translate of the set $A$.    We wish to obtain bounds for $\P(A)$ and $\P(A \cap A(\tau))$.   

 \begin{proposition}
\label{prop: two point couple}
Let $|\tau| \le 1$ be a real number, and let $0\le m\le K-3$ denote the largest integer in this range with $|\tau| \le (\log T)^{-m/K}$.  
Then, for any choice of parameters  $0<x_j\leq \log\log T$ (with $1\le j \le K-3$), we have 
\begin{equation}\label{eqn: two point mesoscopic}
\mathbb P ( A\cap A(\tau) )
\ll \exp\Big( - \sum_{j=1}^{m}\frac{x_j^2}{2s_j^2} - \sum_{j=m+1}^{K-3} \frac{x_j^2}{s_j^2} \Big).
\end{equation}
\end{proposition}

\begin{proof}  For brevity, we write $P_j=P_j(t)$, $P'_j=P_j(t+\tau)$,  $B'=B(\tau)$, and $A'=A(\tau)$.  We shall bound 
$\P(A \cap B \cap A'\cap B')$, and then the bound of the proposition will follow since the complements of the sets $B$ and $B'$ 
have very small measure, by Lemma \ref{lem: barrier}.  

For any choice of parameters $\beta_j >0$ (for $1\le j\le K-3$), for $t$ in the set $A\cap A'$ we have 
$$ 
\sum_{j=1}^{K-3} \beta_j (P_j +P_j^{\prime}) \ge 2 \sum_{j=1}^{K-3} \beta_j x_j. 
$$ 
Therefore 
$$
\P(A \cap B \cap A'\cap B') \le \mathbb{E} \Big[  \exp \Big ( \sum_{j=1}^{K-3} \beta_j ( P_j + P'_j ) \Big)
		  {\bf 1}_{B \cap B'} \Big] \exp\Big( - 2 \sum_{j=1}^{K-3} \beta_j  x_j \Big). 
$$
Assuming that $\beta_j \le (\log T)^{\frac{1}{16K}}$ for all $j$, from \eqref{eqn: two point exp moment} this is 
$$
 \ll \exp \Big( \frac{1}{2} \sum_{j=1}^{K-3} 2\beta_j^2 (s_j^2 + \rho_j(\tau))- 2 \sum_{j=1}^{K-3} \beta_j x_j \Big).
 $$
If $0\le m\le K-3$  denotes the largest integer with $|\tau| \le 2(\log T)^{-m/K}$ then 
for $1\le j\le m$ we have the trivial bound $\rho_j(\tau) \le s_j^2$, and for $K-3 \ge j \geq m+1$ 
we have by \eqref{5.01} that $\rho_j(\tau) = \OO(1)$.  Therefore our bound above is 
$$
\ll \exp \Big(
\frac{1}{2} \sum_{j=1}^{m} 4\beta_j^2 s_j^2
+ \frac{1}{2} \sum_{j=m+1}^{K-3} (2\beta_j^2 s_j^2 + \OO (\beta_j^2))
- 2 \sum_{j=1}^{K-3} \beta_j x_j \Big).
$$
By setting $\beta_j = x_j/s_j^2$ for $j \ge m+1$ and $\beta_j = x_j/(2s_j^2)$ for $j \le m$ we obtain
\eqref{eqn: two point mesoscopic}.
\end{proof}

The crude bound of Proposition \ref{prop: two point couple} will be sufficient when $|\tau| \le (\log T)^{-\frac 1{2K}}$, 
but when $|\tau|$ is larger (almost of macroscopic size) we will require more  
precise large deviation bounds.   
These can be obtained by doing a change of measure under which the value $x_j$ is typical for $P_j$, and by applying a Berry-Esseen type bound. 
This approach was taken in \cite{ArgBelBou15}.
We use a different approach here by directly inverting the Fourier transform. 
To state the results cleanly, it is convenient to set 
\begin{equation} 
\label{5.10} 
\Psi(x) = \frac{1}{\sqrt{2\pi }}\int_{x}^{\infty} e^{-y^2/2} \rd y, 
\end{equation} 
which is the probability of a standard normal random variable being larger than $x$.  

\begin{proposition}
\label{prop: prob decouple}
For all choices of $0<x_j\leq \log\log T$ (with $1\le j\le K-3$) we have 
\begin{equation}
\label{eqn: prob one point}
\mathbb P(A ) = (1+\oo(1))
\prod_{j=1}^{K-3} \Psi(x_j/s_j).
\end{equation} 
Moreover, if $1 \ge |\tau| \geq (\log T)^{-\frac{1}{2K}}$, then
\begin{equation}
\label{eqn: prob two points}
\mathbb P(A \cap A(\tau) )
= (1+\oo(1))\ \mathbb P(A )\ \mathbb P(A(\tau) )= (1+\oo(1)) \P (A)^2.
\end{equation}
\end{proposition}

\begin{proof}
The proof is based on 
 inverting the Fourier-Laplace transform and using the work in Proposition \ref{prop: fourier}.  We 
 begin with a simple, but useful, contour integral.   Let $x$ be a real number, and $c$ be positive. 
 Then 
 $$ 
 \frac{1}{2\pi i} \int_{c-i\infty}^{c+i\infty} \frac{e^{xw} }{w^2} \rd w = \begin{cases} 
 x &\text{ if  } x\ge 0\\ 
 0 &\text{if  } x \le 0.\\
 \end{cases}
 $$ 
 This may be proved by shifting the contour to the right for $x \le 0$, and to the left (picking up the 
 contribution of the pole at $w=0$) when $x>0$.    Now let $\delta$ be a positive real 
 number.  Applying the identity above twice we find 
  \begin{equation} 
 \label{5.11} 
 \frac{1}{2\pi i} \int_{c-i\infty}^{c+i\infty} e^{xw} \frac{e^{\delta w}-1}{\delta w} \frac{\rd w}{w} = \begin{cases} 
 1 &\text{ if } x\ge 0\\ 
 (\delta+ x)/\delta &\text{ if } -\delta \le x\le 0 \\ 
 0 &\text{ if } x \le {-\delta}. 
 \end{cases} 
 \end{equation}   
 Call the function on the right side above $g_\delta(x)$, which plainly approximates the indicator function of the positive reals: 
 ${\bf 1}_{x\ge 0} \le g_\delta(x) \le {\bf 1}_{x+\delta \ge 0}$.

We use the same notation $P_j$, $P_j'$, $A$, $A'$, $B$, $B'$ as in Proposition \ref{prop: two point couple}.
We start with the one-point bound \eqref{eqn: prob one point}.  Since the measure of $B^c$ is negligible, it suffices to evaluate $\mathbb{P}(A \cap B)$.  
We take $\delta = (\log T)^{-\frac 1{64K^2}}$, and from the definition of $g_\delta$  we see that 
\begin{align} 
\label{5.12} 
\P (A\cap B)  &\le {\E} \Big[ \prod_{j} g_\delta (P_j -x_j) {\bf 1}_B \Big] \nonumber \\ 
&= \frac{1}{(2\pi i)^{K-3}} \int_{w_j, \text{Re}(w_j) =\beta_j}  \E \Big[ \exp\Big( \sum_{j} w_j(P_j -x_j)\Big)  {\bf 1}_B\Big]
\prod_{j} \Big( \frac{e^{\delta w_j}-1}{\delta w_j}\Big) \frac{\rd w_j}{w_j},
\end{align}
where we have a $(K-3)$-fold integral with the variables $w_j$ lying on the lines with \text{Re}$(w_j) =\beta_j$ 
with $\beta_j = x_j/s_j^2$.  Note that  $1/\log \log T \le \beta_j  \le (\log \log T)/s_j^2 = \OO(1)$.   

To evaluate the integral above, we draw on our 
work in Proposition \ref{prop: fourier} which will apply when all the $|w_j|$ are bounded by $(\log T)^{\frac 1{16K}}$.  We 
first bound the contribution from terms where some $|w_j|$ is larger than $(\log T)^{\frac 1{16K}}$.   Since 
\begin{align*}
\Big| \E\Big[ \exp \Big (\sum_j w_j (P_j-x_j)\Big) {\bf 1}_B \Big]  \Big | & \le \E \Big[ \exp \Big ( \sum_{j} \beta_j (P_j-x_j) \Big ) \mathbf{1}_{B} \Big]  \\
&\ll \exp\Big( \sum_{j} \Big(\frac{\beta_j^2 s_j^2}{2} - \beta_j x_j\Big) \Big) = 
\exp\Big( -\sum_{j} \frac{x_j^2}{2s_j^2} \Big), 
\end{align*}
such terms contribute (assuming that $|w_1| > (\log T)^{\frac 1{16K}}$, the other cases being similar) 
$$ 
\ll \int_{\substack{ \text{Re}(w_j) =\beta_j \\  |w_1| \ge (\log T)^{\frac{1}{16K}}} }
 \exp\Big( - \sum_{j} \frac{x_j^2 }{2s_j^2}  \Big) \prod_{j} \frac{1}{\delta |w_j|^2} |\rd w_j| 
 \ll  \exp\Big( - \sum_{j} \frac{x_j^2 }{2s_j^2}  \Big)  \prod_{j} \frac{1}{\delta \beta_j}  (\log T)^{-\frac{1}{16K}}. 
$$ 
Recalling that $\delta = (\log T)^{-\frac{1}{64K^2}}$, a small calculation using $\Psi(x) \gg e^{-x^2/2}/(1+x)$ for all $x\ge 0$ 
shows that the above is 
\begin{equation} 
\label{5.13} 
\ll (\log T)^{-\frac{1}{32K}} \prod_{j} \Psi(x_j/s_j). 
\end{equation}

Now we turn to the portion of the integral in \eqref{5.12} where all the variables $w_j$ are bounded in size by $(\log T)^{\frac{1}{16K}}$.  
Here we use \eqref{eqn: one point exp moment}, and obtain 
$$ 
\frac{1}{(2\pi i)^{K-3}} \int_{\substack{ w_j, \text{Re}(w_j) =\beta_j \\ |w_j| \le (\log T)^{\frac 1{16K}}}} \Big(\exp\Big(\sum_j \Big(\frac{w_j^2 s_j^2 }{2} - w_j x_j \Big) 
\Big) + \OO (\exp(-(\log T)^{\frac{1}{4K}})) \Big)  \prod_{j} \Big( \frac{e^{\delta w_j}-1}{\delta w_j}\Big) \frac{\rd w_j}{w_j}.
$$  
The error term above contributes 
$\ll \exp(-(\log T)^{\frac{1}{4K}}) \prod_{j} 1/(\delta \beta_j) $, which is much smaller than \eqref{5.13}.  
In the main term above we extend the integrals to all ranges of $w_j$, incurring an error bounded again by \eqref{5.13}.  
We are left to handle 
\begin{equation} 
\label{5.14}
\frac{1}{(2\pi i)^{K-3}} \int_{ w_j, \text{Re}(w_j) =\beta_j} \exp\Big(\sum_j \Big(\frac{w_j^2 s_j^2 }{2} - w_j x_j \Big)  \Big)  \prod_{j} \Big( \frac{e^{\delta w_j}-1}{\delta w_j}\Big) \frac{\rd w_j}{w_j}.
\end{equation} 
 If $X_j$ denotes a Gaussian random variable with mean $0$ and variance $s_j^2$, chosen independently for different $j$, then this integral equals 
 \begin{align*} 
 \frac{1}{(2\pi i)^{K-3}} & \int_{ w_j, \text{Re}(w_j)  =\beta_j}  \E\Big[ \exp\Big( \sum_{j} w_j (X_j -x_j)\Big)\Big ]  \prod_{j} \Big( \frac{e^{\delta w_j}-1}{\delta w_j}\Big) \frac{\rd w_j}{w_j} 
 \\
 &= \E \Big[ \prod_j g_\delta(X_j -x_j) \Big ] \le \prod_j \Psi\Big( \frac{x_j}{s_j} - \delta\Big) = (1+\OO(\delta^{1/2})) \prod_{j} \Psi\Big(\frac{x_j}{s_j}\Big).
 \end{align*}
 
 Putting together our analysis, we conclude that $\P (A\cap B) \le (1+\oo(1)) \prod_{j} \Psi(x_j/s_j)$, obtaining the upper bound implicit in 
 \eqref{eqn: prob one point}.   The corresponding lower bound follows similarly starting with $\P (A\cap B) \ge \E[ \prod_j g_{\delta}(P_j-x_j -\delta) {\bf 1}_B]$. 
 
 The proof of \eqref{eqn: prob two points} is similar.   Here we start with 
 \begin{align} 
\label{5.15} 
&\P (A\cap A' \cap B\cap B') \le {\E} \Big[ \prod_{j} g_\delta (P_j -x_j) g_\delta(P_j' -x_j) {\bf 1}_{B\cap B'} \Big] \nonumber \\ 
&= \frac{1}{(2\pi i)^{2(K-3)}} \int_{\substack{ w_j , w_j' \\  \text{Re}(w_j) = \text{Re}(w_j')  =\beta_j}} 
  \E \Big[ \exp\Big( \sum_{j} ( w_j(P_j -x_j) +w_j' (P_j'-x_j) ) \Big)  {\bf 1}_{B\cap B'}\Big] \nonumber \\ 
&\hskip 1.5 in \times \prod_{j} \Big( \frac{e^{\delta w_j}-1}{\delta w_j}\Big) \Big( \frac{e^{\delta w_j'}-1}{\delta w_j'}\Big)  \frac{\rd w_j}{w_j} \frac{\rd w_j'}{w_j'}.  
\end{align}
We proceed as before, bounding the tails of the integrals where some $w_j$ or $w_j'$ has size $> (\log T)^{\frac{1}{16K}}$ as we did in 
\eqref{5.13}.   For the remaining integrals with $|w_j|$ and $|w_j'|  \le  (\log T)^{\frac{1}{16K}}$ we use \eqref{eqn: two point exp moment} of 
Proposition \ref{prop: fourier}.  After estimating the error terms arising here, and extending the integrals over $w_j$ and $w_j'$ (exactly as before) we 
arrive, in place of \eqref{5.14}, at 
\begin{align} 
\label{5.16} 
\frac{1}{(2\pi i)^{2(K-3)}} \int_{\substack{ w_j, w_j' \\ 
 \text{Re}(w_j) = \text{Re}(w_j') =\beta_j} } &
 \exp\Big(\sum_j \Big(\frac{w_j^2 s_j^2 }{2} - w_j x_j  + \frac{(w_j's_j)^2}{2} - w_j'x_j  +\rho_j(\tau) w_j w_j'\Big)  \Big) \nonumber\\ 
 &\times \prod_{j} \Big( \frac{e^{\delta w_j}-1}{\delta w_j}\Big) \Big( \frac{e^{\delta w_j'}-1}{\delta w_j'}\Big) \frac{\rd w_j}{w_j} \frac{\rd w_j'}{w_j'}.
\end{align} 
 Since $|\tau| \ge (\log T)^{-\frac{1}{2K}}$, from \eqref{5.02} we have $\rho_j(\tau) = \OO ( (\log T)^{-\frac 1{2K}})$ for all $j$ and 
 therefore the cross terms $\exp(\rho_j(\tau) w_j w_j')$ appearing in \eqref{5.16} make a negligible contribution.   We are then left 
 with essentially two copies of the integrals in \eqref{5.14}, enabling us to conclude that 
 $$ 
 \P   (A\cap A' \cap B\cap B') \le (1+\oo(1)) \prod_{j} \Psi(x_j/s_j)^2. 
 $$ 
 As before, we can obtain the corresponding lower bound as well, completing the proof of \eqref{eqn: prob two points}.  
\end{proof}

\subsection{Proof of Proposition \ref{prop: lower bound dirichlet}}
\label{sect: proof Dirichlet}  Divide the interval $[-1/4,1/4]$ into $\lfloor \log T \rfloor$ equally spaced points $\tau_\ell$  (with $1\le \ell \le 
\lfloor \log T\rfloor$).   Take $x_j = (\lambda \log \log T)/K$ in the definition of the event $A(\tau)$, so that Proposition \ref{prop: lower bound dirichlet} follows if we can establish that 
\begin{equation} 
\label{5.20} 
\P \Big( \bigcup_{\ell} A(\tau_\ell) \Big)  = 1+ \oo (1). 
\end{equation} 
Recall that $A(\tau_\ell)$ is essentially a translate of the set $A$, and so by \eqref{5.10} and \eqref{eqn: prob one point} we 
have 
\begin{equation} 
\label{5.21} 
\P(A(\tau_\ell))  = \P(A) + \OO(1/T)  = (1+\oo(1)) \prod_{j=1}^{K-3} \Psi\Big(\frac {x_j}{s_j}\Big). 
\end{equation} 
Since $s_j^2 = (\log \log T)/(2K) + \OO(1)$, from our choice of $x_j$ and since $\Psi(x) \gg e^{-x^2/2}/x$ for $x\ge 1$, we obtain 
\begin{align} 
\label{5.22}  
\P(A(\tau_\ell))&= (1+\oo(1))  
 \prod_{j=1}^{K-3} 
 \Psi\Big( \lambda \frac{\sqrt{2\log \log T}}{\sqrt{K}}\Big) 
 \gg \prod_{j=1}^{K-3} \frac{(\log T)^{-\lambda^2/K}}{\sqrt{\log \log T}} \nonumber \\
 &= (\log T)^{-\lambda^2 (1-3/K)} (\log \log T)^{-(K-3)/2}. 
\end{align}

The Cauchy-Schwarz inequality gives 
\begin{equation} 
\label{5.23}
\Big( \E \Big[ \sum_{\ell} {\bf 1}_{A(\tau_\ell)} \Big] \Big)^2 = 
\Big( \E \Big[ {\bf 1}_{\cup_{\ell} A(\tau_\ell)}  \sum_{\ell} {\bf 1}_{A(\tau_\ell)} \Big] \Big)^2  
\le \P\Big( \bigcup_{\ell} A(\tau_\ell) \Big) \E \Big[ \Big( \sum_{\ell} {\bf 1}_{A(\tau_\ell)}\Big)^2\Big];  
\end{equation} 
this may be viewed as a special case of the Paley-Zygmund inequality.  Note that, by \eqref{5.21} and \eqref{5.22},   
\begin{equation} 
\label{5.24} 
\Big(\E \Big[ \sum_{\ell} {\bf 1}_{A(\tau_\ell)} \Big] \Big)^2 = \Big( \sum_{\ell} \P (A(\tau_\ell)) \Big)^2  
= \Big( (1+\oo (1)) \lfloor \log T\rfloor \P(A)\Big)^2 \gg (\log T)^{2 (1-\lambda^2 +3 \lambda^2/K) -\epsilon},
\end{equation} 
for $\epsilon>0$.
To establish \eqref{5.20} we now establish an upper bound for the second factor on the right side of \eqref{5.23}.  

Expanding out, we have 
\begin{align} 
\label{5.30}
\E \Big[ \Big( \sum_{\ell} {\bf 1}_{A(\tau_\ell)}\Big)^2\Big] &= \sum_{\ell, \ell'} \P \big( A(\tau_{\ell} ) \cap A(\tau_{\ell'}) \big) \nonumber \\ 
&=\Big( \sum_{|\tau_{\ell} - \tau_{\ell'}| \ge (\log T)^{-1/(2K)} } + \sum_{|\tau_{\ell} - \tau_{\ell'}| \le (\log T)^{-1/(2K)} } \Big) \P \big( A(\tau_{\ell} ) \cap A(\tau_{\ell'}) \big).
\end{align} 
The first term accounts for the typical pair of points $\tau_{\ell}$, $\tau_{\ell'}$, and by \eqref{eqn: prob two points} we have (for 
such a pair) 
$$ 
 \P \big( A(\tau_{\ell} ) \cap A(\tau_{\ell'}) \big) =  \P \big( A \cap A(\tau_\ell - \tau_{\ell'}) \big) +\OO(1/T)  = (1+\oo(1)) \P(A(\tau_\ell)) \P(A(\tau_\ell')). 
$$
Therefore 
\begin{equation} 
\label{5.31} 
 \sum_{|\tau_{\ell} - \tau_{\ell'}| \ge (\log T)^{-1/(2K)} }   \P \big( A(\tau_{\ell} ) \cap A(\tau_{\ell'}) \big) \le (1+\oo(1)) \Big( \sum_{\ell} 
 \P(A(\tau_\ell))\Big)^2. 
 \end{equation} 

We now bound the second term in \eqref{5.30}, using Proposition \ref{prop: two point couple} to show that its contribution is negligible.  
Let $m$ denote the largest integer in $[0,K-3]$ with $|\tau_{\ell}-\tau_{\ell'}| \le (\log T)^{-m/K}$.   Then Proposition \ref{prop: two point couple} gives (since $x_j^2/s_j^2 =2 \lambda^2 (\log \log T)/K$) 
$$ 
\P \big( A(\tau_{\ell} ) \cap A(\tau_{\ell'}) \big)  \ll  (\log T)^{- \lambda^2 (m/K + 2(K-3-m)/K)}. 
$$
In the range $1\le m\le K-3$, the number of pairs ($\tau_{\ell}$, $\tau_{\ell'}$) is $\ll (\log T)^{2-m/K}$, while in the case $m=0$ (since 
we are considering the case $|\tau_{\ell} -\tau_{\ell'}| \le (\log T)^{-\frac{1}{2K}}$) the number of pairs is $\ll (\log T)^{2-\frac{1}{2K}}$.  
It follows that 
\begin{align*}
 \sum_{|\tau_{\ell} - \tau_{\ell'}| \le (\log T)^{-1/(2K)} }   \P \big( A(\tau_{\ell} ) \cap A(\tau_{\ell'}) \big) 
 &\ll (\log T)^{2-\frac{1}{2K}} (\log T)^{-2\lambda^2 (1-3/K)} \\
 &+ \sum_{m=1}^{K-3} (\log T)^{2-m/K} (\log T)^{-2\lambda^2 (1-3/K) + m\lambda^2/K}. 
 \end{align*}
 Since $\lambda <1$, using \eqref{5.24} (with $\epsilon$ there sufficiently small) we conclude that 
 \begin{equation} 
 \label{5.32} 
  \sum_{|\tau_{\ell} - \tau_{\ell'}| \le (\log T)^{-1/(2K)} }   \P \big( A(\tau_{\ell} ) \cap A(\tau_{\ell'}) \big)  = \oo\Big( \Big( \sum_{\ell} 
 \P(A(\tau_\ell))\Big)^2 \Big). 
 \end{equation}  
 
 From \eqref{5.30}, \eqref{5.31}, and \eqref{5.32} we conclude that 
 $$ 
 \E \Big[ \Big( \sum_{\ell} {\bf 1}_{A(\tau_\ell)}\Big)^2\Big] \le (1+\oo(1))  \Big( \sum_{\ell} 
 \P(A(\tau_\ell))\Big)^2.
 $$ 
 Inserting this upper bound in \eqref{5.23}, and using \eqref{5.24}, we deduce the bound \eqref{5.20}, which 
 completes the proof of our proposition.


\begin{bibdiv}
\begin{biblist}


\bib{Arg16}{article}{
    AUTHOR = {Arguin, L.-P.},
     TITLE = {Extrema of Log-correlated Random Variables: Principles and Examples},
   JOURNAL = {in Advances in Disordered Systems, Random Processes and Some Applications, Cambridge Univ. Press, Cambridge},
      YEAR = {2016},
        PAGES = {166--204}
}

\bib{ArgBelBou15}{article}{
    AUTHOR = {Arguin, L.-P.},
    AUTHOR = {Belius ,D.},
    AUTHOR = {Bourgade, P.},
     TITLE = {Maximum of the characteristic polynomial of random unitary matrices},
   JOURNAL = {{\it Comm. Math. Phys.}},
    VOLUME = {349},
      YEAR = {2017},
      PAGES = {703--751}
}


\bib{ArgBelHar15}{article}{
    AUTHOR = {Arguin, L.-P.},
    AUTHOR = {Belius ,D.},
    AUTHOR = {Harper, A.J.},
     TITLE = {Maxima of a randomized Riemann zeta function, and branching random walks},
     JOURNAL = {Ann. Appl. Probab.},
  FJOURNAL = {The Annals of Applied Probability},
    VOLUME = {27},
      YEAR = {2017},
    NUMBER = {1},
     PAGES = {178--215},
      ISSN = {1050-5164},
   MRCLASS = {60G70 (11M06)},
  MRNUMBER = {3619786},
   
}



\bib{BalRam77}{article}{
    AUTHOR = {Balasubramanian, R.},
    AUTHOR = {Ramachandra, K.},
     TITLE = {On the frequency of Titchmarsh'��s phenomenon for $\zeta(s)$-III},
   JOURNAL = {Proc. Indian Acad. Sci.},
    VOLUME = {86},
      YEAR = {1977},
     PAGES = {341--351}
}




\bib{BhaRao76}{book}{
   author={Bhattacharya, R. N.},
   author={Ranga Rao, R.},
   title={Normal approximation and asymptotic expansions},
   note={Wiley Series in Probability and Mathematical Statistics},
   publisher={John Wiley \& Sons, New York-London-Sydney},
   date={1976}
}

\bib{BoGh}{article}{
  author = {Bober, J. W.},
  author = {Hiary, G. A.},
  title = {New Computations of the Riemann Zeta Function on the Critical Line},
  journal = {Exp. Math.},
  year  = {2016},
  publisher = {Taylor & Francis},
  eprint = {
    http://dx.doi.org/10.1080/10586458.2016.1233083
  }
}

\bib{BonSei15}{article}{
    AUTHOR = {Bondarenko, A.},
   AUTHOR = {Seip, K.},
     TITLE = {Large GCD sums and extreme values of the Riemann zeta function},
  journal={Duke Math. J.},
   volume={166},
   date={2017},
   number={9},
   pages={1685--1701},
  }



\bib{Bou16}{article}{
    AUTHOR = {Bourgain, J.},
     TITLE = {Decoupling, exponential sums and the Riemann zeta function},
    journal={J. Amer. Math. Soc.},
   volume={30},
   date={2017},
   number={1},
   pages={205--224},
   }


\bib{Bou10}{article}{
    AUTHOR = {Bourgade, P.},
     TITLE = {Mesoscopic fluctuations of the zeta zeros},
   JOURNAL = {Probab. Theory Related Fields},
    VOLUME = {148},
      YEAR = {2010},
    NUMBER = {3-4},
     PAGES = {479--500}
}



\bib{Bra78}{article}{
   author={Bramson, M. D.},
   title={Maximal displacement of branching Brownian motion},
   journal={Comm. Pure Appl. Math.},
   volume={31},
   date={1978},
   number={5},
   pages={531--581}
}



\bib{ChaSou11}{article}{
    AUTHOR = {Chandee, V.},
    AUTHOR = {Soundararajan,K.},
     TITLE = {Bounding $|\zeta(1/2+\ii t)|$ on the Riemann hypothesis},
   JOURNAL = {Bull. London Math. Soc.},
      YEAR = {2011},
}



\bib{ChaMadNaj16}{article}{
    AUTHOR = {Chhaibi, R.},
    AUTHOR = {Madaule,T.},
    AUTHOR = {Najnudel, J.},
     TITLE = {On the maximum of the C$\beta$E field},
   JOURNAL = {Preprint arXiv:1607.00243},
      YEAR = {2016},
}

\bib{Con88}{article}{
    AUTHOR = {Conrey, J. B.},
     TITLE = {The fourth moment of derivatives of the {R}iemann
              zeta-function},
   JOURNAL = {Quart. J. Math. Oxford Ser. (2)},
  FJOURNAL = {The Quarterly Journal of Mathematics. Oxford. Second Series},
    VOLUME = {39},
      YEAR = {1988},
    NUMBER = {153},
     PAGES = {21--36},

}

\bib{FarGonHug07}{article}{
  title = {The maximum size of L-functions},
  author = {Farmer, D.},
  author = {Gonek, S.M.},
  author = {Hughes, C.},
  journal = {J. Reine Angew. Math.},
  volume = {609},
  pages = {215--236},
  year = {2007},
}




\bib{FyoBou08}{article}{
    AUTHOR = {Fyodorov, Y. V.},
    AUTHOR = {Bouchaud, J.-P.},
     TITLE = {Freezing and extreme-value statistics in a random energy model
              with logarithmically correlated potential},
   JOURNAL = {J. Phys. A},
    VOLUME = {41},
      YEAR = {2008},
    NUMBER = {37},
     PAGES = {372001, 12 pp.}
}


\bib{FyoHiaKea12}{article}{
  title = {Freezing Transition, Characteristic Polynomials of Random Matrices, and the {Riemann} Zeta Function},
  author = {Fyodorov, Y. V.},
  author = {Hiary, G. A.},
  author = {Keating, J. P.},
  journal = {Phys. Rev. Lett.},
  volume = {108},
  pages = {170601, 5pp.},
  year = {2012},
  publisher = {American Physical Society}
}



\bib{FyoKea14}{article}{
    AUTHOR = {Fyodorov, Y. V.},
    AUTHOR = {Keating, J. P.},
     TITLE = {Freezing transitions and extreme values: random matrix theory,
              and disordered landscapes},
   JOURNAL = {Phil. Trans. R. Soc. A},
    VOLUME = {372},
      YEAR = {2014},
     PAGES = {20120503, 32 pp. }
}



\bib{FyoLedRos09}{article}{
    AUTHOR = {Fyodorov, Y. V.},
    AUTHOR = {Le Doussal, P.},
author={Rosso, A.},
     TITLE = {Statistical mechanics of logarithmic {REM}: duality, freezing
              and extreme value statistics of {$1/f$} noises generated by
              {G}aussian free fields},
   JOURNAL = {J. Stat. Mech. Theory Exp.},
      YEAR = {2009},
    NUMBER = {10},
     PAGES = {P10005, 32 pp.}
}




\bib{FyoLedRos12}{article}{
   author={Fyodorov, Y. V.},
   author={Le Doussal, P.},
   author={Rosso, A.},
   title={Counting function fluctuations and extreme value threshold in
   multifractal patterns: the case study of an ideal $1/f$ noise},
   journal={J. Stat. Phys.},
   volume={149},
   date={2012},
   number={5},
   pages={898--920},
}



\bib{FyoSim2015}{article}{
    AUTHOR = {Fyodorov, Y. V.},
    AUTHOR = {Simm, N. J.},
     TITLE = {On the distribution of maximum value of the characteristic polynomial of GUE random matrices},
    journal={Nonlinearity},
   volume={29},
   date={2016},
   number={9},
   pages={2837--2855},
}

\bib{Har13}{article}{
    AUTHOR = {Harper, A. J.},
     TITLE = {A note on the maximum of the Riemann zeta function, and log-correlated random variables},
   JOURNAL = {Preprint arXiv:1304.0677},
      YEAR = {2013},
}


\bib{Ing28}{article}{
   author={Ingham, A. E.},
   title={Mean-Value Theorems in the Theory of the Riemann Zeta-Function},
   journal={Proc. London Math. Soc. (2)},
   volume={27},
   date={1928},
   pages={273--300}
   }

\bib{Ivi03}{book}{
    AUTHOR = {Ivi\'c, A.},
     TITLE = {The {R}iemann zeta-function},
      NOTE = {Theory and applications},
 PUBLISHER = {Dover Publications, Inc., Mineola, NY},
      YEAR = {2003},
     PAGES = {xxii+517},
      ISBN = {0-486-42813-3},
   MRCLASS = {11M06 (11M26)},
  MRNUMBER = {1994094},
}



\bib{KeaSna00}{article}{
   author={Keating, J.},
   author={Snaith, N.},
   title={Random Matrix Theory and $\zeta(1/2+\ii t)$},
   journal={Comm. Math. Phys.},
   number={1},
   pages={57--89},
   date={2000}
}





\bib{Kis15}{collection}{
    AUTHOR = {Kistler, N.},
     TITLE = {Derrida's Random Energy Models},
 BOOKTITLE = {Correlated Random Systems: Five Different Methods},
    SERIES = {Lecture Notes in Math.},
    VOLUME = {2143},
     PAGES = {71--120},
 PUBLISHER = {Springer, Berlin},
      YEAR = {2015},
}




\bib{LamPaq16}{article}{
author = {Lambert, G.},
author = {Paquette, E.},
journal = {preprint, arXiv:1611.08885},
title = {The law of large numbers for the maximum of almost Gaussian log-correlated fields coming from random matrices},
year = {2016}
}



\bib{Lin08}{article}{
author = {Lindel\"{o}f, E.},
journal = {Bull. Sci. Math.},
pages = {341--356},
title = {Quelques remarques sur la croissance de la fonction $\zeta(s)$},
volume = {32},
year = {1908},
}


\bib{Lit24}{article}{
author = {Littlewood, J.E.},
journal = {Proc. Camb. Philos. Soc.},
pages = {295--318},
title = {On the zeros of the Riemann zeta-function},
volume = {22},
year = {1924}
}


\bib{Mon73}{article}{
    AUTHOR = {Montgomery, H.},
     TITLE = {The pair correlation of zeros of the zeta function},
   conference={
      title={Analytic number theory},
      address={Proc. Sympos. Pure Math., Vol. XXIV, St. Louis Univ., St.
      Louis, Mo.},
      date={1972},
   },
   book={
      publisher={Amer. Math. Soc., Providence, R.I.},
   },
   date={1973},
   pages={181--193}
}



\bib{Mon77}{article}{
    AUTHOR = {Montgomery, H.},
     TITLE = {Extreme values of the Riemann zeta function},
   journal={Comment. Math. Helv.},
   volume={52},
   pages={511--518},      
   YEAR = {1977},
}

\bib{MoVa07}{book}{
    AUTHOR = {Montgomery, H.},
      AUTHOR = {Vaughan, R.},
     TITLE = {Multiplicative number theory. {I}. {C}lassical theory},
    SERIES = {Cambridge Studies in Advanced Mathematics},
    VOLUME = {97},
 PUBLISHER = {Cambridge University Press, Cambridge},
      YEAR = {2007},
     PAGES = {xviii+552},

}

\bib{Naj16}{article}{
author = {Najnudel, J.},
journal = {Preprint arXiv:1611.05562},
title = {On the extreme values of the Riemann zeta function on random intervals of the critical line},
year = {2016}
}


\bib{PaqZei17}{article}{
author = {Paquette, E.},
author = {Zeitouni, O.},
journal = {International Mathematics Research Notices},
title = {The maximum of the CUE field},
year = {2017},
pages={1--92},
}


\bib{Rad11}{article}{
author = {Radziwill, M.},
journal = {Preprint arxiv:1108.5092},
title = {Large deviations in Selberg's central limit theorem},
year = {2011}
}



\bib{RadSou15}{article}{
author = {Radziwill, M.},
author = {Soundararajan, K.},
journal = {Preprint arxiv:1509.06827},
title = {Selberg's central limit theorem for $\log|\zeta(\frac{1}{2}+\ii t)|$},
year = {2015}
}



\bib{RamSan93}{article}{
author = {Ramachandra, K.},
author = {Sankaranarayanan, A.},
journal = {J. Number Theory},
pages = {281--291},
title = {On some theorems of Littlewood and Selberg},
volume = {44},
year = {1993}
}



\bib{SakWeb16}{article}{
author = {Saksman, E.},
author = {Webb, C.},
journal = {preprint, arXiv:1609.00027},
title = {The Riemann zeta function and Gaussian multiplicative chaos: statistics on the critical line},
year = {2016}
}
 
 
\bib{Sel46}{article}{
author = {Selberg, A.},
journal = {Arch. Math. Naturvid.},
pages = {315--392},
publisher = {Arch. Math. Naturvid.},
title = {Contributions to the theory of the Riemann zeta-function},
volume = {48},
pages={89--155},
year = {1946},
}


\bib{Sou08}{article}{
author = {Soundararajan, K.},
title={Extreme values of zeta and L-functions},
journal = {Mathematische Annalen},
year = {2008},
volume={342},
issue={2},
pages={467--486}
}

\bib{Sou09}{article}{
author = {Soundararajan, K.},
title={Moments of the Riemann zeta-function},
journal = {Annals of Math.},
volume={170},
pages={981--993},
year = {2009}
}



\bib{Tit86}{article}{
author = {Titchmarsh, E.C.},
journal = {Second Edition, Oxford Univ. Press, New
York},
title = {The theory of the Riemann zeta-function},
year = {1986}
}


\end{biblist}
\end{bibdiv}
\end{document}